\documentclass[11pt]{article}

\usepackage[nocompress]{cite}
\usepackage{amsfonts}
\usepackage[fleqn,reqno]{amsmath}
\usepackage{amssymb}
\usepackage[titletoc]{appendix}
\usepackage{comment}
\usepackage{enumitem}
\usepackage{filecontents}
\usepackage[top=1.2in,bottom=1.2in,left=1in, right=1in]{geometry}
\usepackage{graphics}
\usepackage{lineno}
\usepackage{pgfplots}
\usepackage{tikz}
\usepackage{todonotes}
\usepackage{placeins}
\usepackage{bm}
\usepackage[linesnumbered,ruled]{algorithm2e}
\usetikzlibrary{arrows}

\usepackage[pagebackref=false,bookmarks=false]{hyperref} 

\hypersetup{
  bookmarksnumbered=true,
  bookmarksopen=false,
  hypertexnames=false,      
  breaklinks=true,          
  unicode=false,
  pdffitwindow=true,        
  pdfnewwindow=true,        
  colorlinks=true,         
  linkcolor=blue,
  anchorcolor=red,
  citecolor=blue,
  filecolor=magenta,
  urlcolor=blue,
  pdfstartview = FitH,
  pdfkeywords = {},
  pdfcreator = {LaTeX with hyperref package}
}

\providecommand{\e}[1]{\ensuremath{\times 10^{#1}}}

\newcommand{\bd}{{\partial}}

\newcommand{\cc}{{\mathbf{c}}}
\newcommand{\dd}{{\mathbf{d}}}
\newcommand{\DD}{{\mathcal{D}}}
\newcommand{\eeta}{{\boldsymbol\eta}}

\newcommand{\FF}{{\mathbf{F}}}
\newcommand{\grad}{{\nabla}}

\newcommand{\nn}{{\mathbf{n}}}

\newcommand{\rr}{{\mathbf{r}}}
\renewcommand{\Re}{{\mathrm{Re}}}

\newcommand{\ssigma}{{\boldsymbol\sigma}}
\renewcommand{\tt}{{\mathbf{t}}}
\newcommand{\uu}{{\mathbf{u}}}
\newcommand{\UU}{{\mathbf{U}}}

\newcommand{\xx}{{\mathbf{x}}}

\newcommand{\yy}{{\mathbf{y}}}

\def\gap{\hspace*{.2in}}

\begin{document}

\title{Stable and contact-free time stepping for dense rigid particle
suspensions}

\author{Lukas Bystricky, Sachin Shanbhag, Bryan D.~Quaife}
\date{Department of Scientific Computing, Florida State University}


\maketitle

\begin{abstract} 
We consider suspensions of rigid bodies in a two-dimensional viscous
fluid. Even with high-fidelity numerical methods, unphysical contact
between particles occurs because of spatial and temporal discretization
errors.  We apply the method of Lu et al.~[{\em Journal of Computational
Physics}, {\bf 347}:160--182, 2017] where overlap is avoided by
imposing a minimum separation distance.  In its original form, the
method discretizes interactions between different particles explicitly.
Therefore, to avoid stiffness, a large minimum separation distance is
used.  In this paper, we extend the method of Lu et al.~by treating all
interactions implicitly.  This new time stepping method is able to
simulate dense suspensions with large time step sizes and a small
minimum separation distance.  The method is tested on various unbounded
and bounded flows, and rheological properties of the resulting
suspensions are computed.

\end{abstract}

\section{Introduction\label{s:intro}}
Dispersions of particulate rods or fibers are used in composite
materials to tune mechanical, thermal, and electrical properties.
Typically, these materials are processed in the melt or liquid
suspension state via operations like injection molding, extrusion, or
casting. It is important to model fiber suspensions for two reasons: (i)
the distribution and orientation of the fibers, which determines the
properties of the composite material, are governed by the flow history
during processing, and (ii) the rheological properties of the
suspension, which influence the flow behavior, in turn, depend on the
size, shape, distribution, and orientation of the
fibers~\cite{larsoncf}.

The theory of rigid fibers in flowing fluids was pioneered by
Jeffery~\cite{Jeffery1922} who analyzed the motion of a single
spheroidal particle sheared in a Newtonian solvent. At a given shear
rate $\dot{\gamma}$, he observed that fibers of length $\ell$ and
diameter $d$ underwent periodic motion with a period $(\pi/\dot{\gamma})
(\lambda + 1/\lambda)$, where $\lambda = \ell/d$ is the aspect ratio.
The period increases with $\lambda$, and when $\lambda \gg 1$, a
particle exhibiting a ``Jeffery's orbit" stays aligned with the flow
direction most of the time, before abruptly spinning through a
half-revolution. In the dilute regime (number of rods/unit volume $\nu <
1/\ell^3$), trajectories of elongated fibers of different shapes, such
as cylinders, can be quantitatively described via Jeffery's orbits once
corrections are made for particle shape~\cite{Bretherton1962}.

As $\nu$ increases, the interactions between fibers become significant.
Batchelor extended Jeffery's theory for multiple particles, by relating
the average stress tensor, ${\bm \sigma}$, to the distribution of fiber
orientation $\mathbf{p}$, and the deformation tensor $\mathbf{D} =
(\nabla \mathbf{u} + \nabla \mathbf{u}^\intercal)/2$. Assuming purely
hydrodynamic interactions between fibers, and a slender body
approximation ($\lambda \gg 1$)~\cite{Batchelor1970, Batchelor1970a,
Doi1978, Dinh1984, Shaqfeh1990},
\begin{align}
  {\bm \sigma} = 2 \mu \mathbf{D} + \nu \zeta 
    \langle \mathbf{p p p p} \rangle : \mathbf{D},
\label{eqn:batchelor}
\end{align}
where $\mu$ is the solvent viscosity, and $\zeta$ is a drag
coefficient~\cite{Batchelor1971} that depends on the size and
concentration of the particles, and the solvent viscosity. The ensemble
average $\langle \cdot \rangle = \int \cdot \,
\psi(\mathbf{p})\,d\mathbf{p}$ represents a weighted average over the
probability distribution of fiber orientations $\psi(\mathbf{p})$. 

In computer simulations, the fiber orientation distribution is modeled
implicitly or explicitly. In the \emph{implicit} approach, individual
fibers are not explicitly represented; instead it relies on averages of
second- and fourth-order fiber orientation tensors, $\langle \mathbf{p
p} \rangle$ and $\langle \mathbf{p p p p} \rangle$. Fluid flow equations
(Stokes or Navier-Stokes) are coupled with evolution equations for the
fiber orientation tensors. In order to solve the resulting equations,
fiber interaction models and closure approximations have to be specified
externally~\cite{Advani1987, Advani1990, Ferec2014, Perez2017}. This is
in contrast to direct numerical simulations where individual fibers are
\emph{explicitly} represented. Typically, fibers are modeled as
prolate ellipsoids~\cite{Ausias2006}, a set of connected
beads~\cite{Yamamoto1996, Joung2001}, rods~\cite{Schmid2000,
Lindstroem2007}, or a slender body ($\ell \gg d$)~\cite{Fan1998,
Rahnama1995, tor-she2004, tor-gus2006, gus-tor2009}, with suitable
first-order corrections to account for finite width. Over the years, in
addition to long-range hydrodynamic interaction, these models have been
supplemented with detailed physics including short-range lubrication,
mechanical contact, and frictional forces~\cite{Sundararajakumar1997,
Lindstroem2008}.

In the semi-dilute regime, $1/\ell^3 \ll \nu \ll 1/d\ell^2$, fiber
rotation is hindered; however, it is found that the statistical
properties are not significantly altered from the dilute
regime~\cite{larsoncf}.  Hydrodynamic interactions between particles
dominate the response, and contacts between fibers are rare. Batchelor's
theory, suitably modified for multibody hydrodynamic
interactions~\cite{Shaqfeh1990, Mackaplow1996}, describes the
empirically observed increase in shear viscosity as a function of $\nu$
reasonably well~\cite{Stover1992, Bibbo1987, Petrich2000}. The
contribution of the fibers to the steady shear viscosity is relatively
modest in non-Brownian suspensions. This is especially true for high
aspect ratio fibers which rotate and align along the flow direction, and
contribute to the viscosity only during the occasional
tumble~\cite{larsoncf}.  Thus, one ought to be careful not to interpret
the success of theory and computer models in predicting the viscosity
change in the semi-dilute regime as validation of the underlying fiber
interaction model. Indeed fiber-fiber interactions are more sensitively
reflected in other viscometric functions such as first normal stress
difference, and distribution of orientations as reflected in, for
example, the dispersion of Jeffery's orbits~\cite{Lindstroem2009}.

Once the concentration increases beyond $\nu \approx 1/d\ell^2$, the
suspension enters the concentrated regime. Here, excluded volume
interactions become important and isotropic packing becomes increasingly
difficult. In this regime, Batchelor's slender body theory and
constitutive relation~\eqref{eqn:batchelor} are no longer valid as
mechanical contacts between fibers start to dominate the response. When
these mechanical interactions are explicitly accounted for, computer
models are able to reproduce a nonzero first normal stress difference
that is observed in experiments~\cite{Sundararajakumar1997, Ausias2006,
Lindstroem2008}. Unlike the dilute and semi-dilute regimes, equation
\eqref{eqn:batchelor} can no longer be used to estimate rheological
properties. Instead, stresses in the suspension have to be computed by
directly summing the forces acting on the fibers~\cite{Ausias2006,
Lindstroem2008}.

In this work, we develop and test tools for two-dimensional direct
numerical simulations of rigid bodies suspended in a viscous fluid.  We
do not make any rigid body assumptions, but rather fully resolve the
fiber shape.  To perform the simulations, we use a boundary integral
equation (BIE) since it resolves the complex geometry by reducing the
set of unknowns to the one-dimensional closed curves that form the fluid
boundary.  Moreover, our BIE fluid solver achieves high-order accuracy.
The governing Stokes equations prohibit contact between particles,
however, because of numerical errors, without additional techniques,
rigid bodies often come into contact or even overlap. Therefore, we
apply a contact algorithm that allows rigid particles to come very close
to one another, but guarantees that contact is avoided without
introducing significant stiffness.  In addition to computing fiber
trajectories, we compute rheological and statistical properties of the
fluid and particles to better understand the dispersion of fibers in
composite materials.

\paragraph{Contributions} Our main contributions are extending the time
stepping strategy introduced for vesicle suspensions~\cite{Quaife2014}
to rigid body suspensions, and analyzing the rheological properties of
the suspensions.  Deformable bodies, such as vesicles, deform as they
approach one another, and this creates a natural minimum separation
distance.  However, for rigid body suspensions, the inability to deform
can force bodies much closer together, and numerical errors can easily
cause unphysical overlap between particles.  To avoid overlap, Lu et
al.~\cite{Lu2017} developed a contact algorithm that guarantees a
minimum separation distance between bodies and use a {\em locally
implicit} time stepping method that only treats inter-body interactions
implicitly.  That is, if $\uu_{ij}$ is the velocity of body $i$ induced
by body $j$, then the time stepping method used is
\begin{align*}
  \frac{\xx_i(t + \Delta t) -  \xx_i(t)}{\Delta t} = 
    \uu_{ii}(t+\Delta t) + \sum_{j \neq i} \uu_{ij}(t),
\end{align*}
where $\xx_i$ is the center of the $i^{\mathrm{th}}$ body.  By treating
the interactions between different bodies explicitly, the minimum
separation distance must be kept sufficiently large to avoid a small
time step restriction due to stiffness---a typical minimum separation
distance is $\mathcal{O}(1)$ arclength spacings.  In line with previous
work of one of the authors~\cite{Quaife2014}, we discretize all
interactions semi-implicitly
\begin{align*}
  \frac{\xx_i(t + \Delta t) -  \xx_i(t)}{\Delta t} = 
    \uu_{ii}(t+\Delta t) + \sum_{j \neq i} \uu_{ij}(t+\Delta t).
\end{align*}
With this modification, we are able to perform simulations with much
smaller and more physical minimum separation distances without
introducing excessive stiffness---a typical minimum separation distance
is $\mathcal{O}(10^{-2})$ arclength spacings.

While maintaining a minimum separation distance is important for stable
simulations, the contact algorithm developed by Lu et al.~does introduce
artificial forces that shifts bodies onto different streamlines, and
this breaks  the reversibility of the Stokes equations.  We examine the
effect of the contact algorithm on the reversibility of the flow.
Finally, we use our new time stepping to examine the rheological
properties of dense rigid body suspensions with small minimum separation
distances.  In particular, we compute the effective shear viscosity of a
suspension of rigid bodies in a Couette device, examine the alignment
angle of elliptical bodies of varying area fraction and aspect ratio,
and compare the results to analytical Jeffery's orbits.

\paragraph{Limitations} The main limitation is that the method is
developed in two dimensions.  By limiting ourselves to two dimensions,
we are able to perform simulations of denser suspensions  than would be
possible in three dimensions.  However, the algorithms we present have
been developed in three dimensions including boundary integral equation
methods and fast summation methods~\cite{cor-gre-rac-vee2017,
kli-tor2014, kli-tor2016}.  The most challenging algorithms to extend to
three dimensions include efficient preconditioners and a suspension
space-time interference volume that integrates a four-dimensional domain
($3$ space dimensions and $1$ time dimension).

\paragraph{Related work} Rather than presenting an exhaustive list of
work related to particulate suspensions in viscous fluids, we focus on
literature related to BIEs and time stepping for rigid body suspensions.
A more complete overview of BIEs for particulate suspensions can be
found in the texts~\cite{Pozrikidis1992, Guazzelli2011, Karrila1991}.
Our work draws heavily from methods developed for simulating
two-dimensional vesicle suspensions~\cite{Quaife2014, Quaife2015,
qua-bir2016, Rahimian2010, Lu2017}.  

We represent the velocity as a completed double-layer potential
representation~\cite{Power1987, Power1993, Karrila1989} that is
discretized with high-order quadrature and solved iteratively with
GMRES~\cite{Saad1986}.  By using a double-layer potential, a second-kind
integral equation needs to be solved.  Upon discretization, the required
number of GMRES iterations is mesh-independent~\cite{Campbell1996}, but
it is geometry-dependent.  Therefore, preconditioners are often applied.
There are a variety of preconditioners available for integral
equations~\cite{cou-pou-dar2017, che2000, qua-cou-dar2018, Quaife2015a,
bra-lub1990, hem-sch1981}, we apply a simple block-diagonal
preconditioner that was successfully used for vesicle
suspensions~\cite{Quaife2014}.

The numerical solution of integral equations requires accurate
quadrature methods for a variety of integrands.  Many of these
integrands are smooth and periodic, and the trapezoid rule is typically
used since it guarantees spectral accuracy~\cite{Trefethan2014}.
However, integrands with large derivatives must be computed when bodies
are in near-contact, and this is a certainty in dense suspensions.  We
apply an interpolation-based quadrature method~\cite{Ying2006,
Quaife2014} since it is efficient and extends to three dimensions, but
other near-singular integration schemes are possible~\cite{Klockner2013,
Barnett2015, Beale2016, Helsing2008, Kropinski1999, Mammoli2006,
Siegel2018}.  The same interpolation-based near-singular integration
scheme is used to compute the pressure and stress, but a combination of
singularity subtraction and odd-even integration~\cite{sid-isr1988,
Quaife2014} is also used to resolve high-order singularities in the
integrands.

The greatest opportunity of acceleration is reducing the cost of the
matrix-vector multiplication required at each GMRES iteration.  We use
the fast multipole method (FMM)~\cite{Greengard1987,Greenbaum1992},
but other fast summation methods, which also extend to three dimensions,
are possible~\cite{bar-hut1986, kli-tor2014}.  As an alternative,
iterations can be entirely avoided by applying a direct solver for
BIEs~\cite{mar-bar-gil-vee2016}, but these solvers would have to be
updated at each time step since the geometry is dynamic.

Once the BIE formulation of the appropriate fluid equations are solved
for the translational and rotational velocities, a time step must be
taken.  We adopt a Lagrangian approach, and since the bodies are rigid,
we only need to track each body's center and inclination angle.
Therefore, for a suspension of $M_p$ bodies, a system of $3M_p$ ordinary
differential equations must be solved---these equations are coupled
through the fluid solver.  Embedded time stepping
methods~\cite{kli-tor2014} work well for dilute suspensions, but can
force the time step to become unreasonably small for moderately dense
suspensions. Artificial repulsion forces~\cite{Flormann2017, Liu2006,
Malhotra2018, Lu2017, Kabacogulu2017} minimize, but do not eliminate,
the chance of a collision. Moreover, these potentials often have sharp
gradients which lead to stiffness and necessitate a small time step
size. Alternatively, a repulsion force based on the concept of
space-time interference volumes (STIVs)~\cite{Harmon2011, Lu2017}
explicitly prevents collisions between particles.  Using current STIV
implementations, the minimum separation distance between bodies cannot
be too small; otherwise, the associated optimization algorithm stalls.
This is a result of treating interactions between different bodies
explicitly.  Therefore, in this work, we extend the STIV contact
algorithm to implicit interactions so that bodies are able to come much
closer---a physical characteristic of dense suspensions of rigid bodies.

In addition to coupling all the bodies implicitly, we further improve
time stepping by allowing for an adaptive time step size.  There are a
variety of adaptive time stepping methods, and they typically either
estimate the local truncation error, and the time step size is adjusted
according to this error~\cite{Quaife2015, Quaife2015a, Sorgentone2018},
or they quantify the computational effort, such as the number of
required time steps, and reduce the time step size when this becomes
large~\cite{Kropinski1999}.  The local truncation error for rigid body
suspensions is expensive because multiple numerical solutions must be
formed.  Therefore, we apply the second option where the time step size
is decreased when the STIV optimization routine requires a large number
of iterations.

\paragraph{Outline of the paper}
In Section~\ref{s:formulation} we describe the physical problem and the
governing equations.  In Section~\ref{s:method}, we describe the
numerical methods used to form numerical solutions.
The results are described in Section~\ref{s:results}, and conclusions
are drawn in Section~\ref{s:conclusions}.

\section{Formulation\label{s:formulation}} 
\begin{figure}[!h]
\begin{center}
\includegraphics{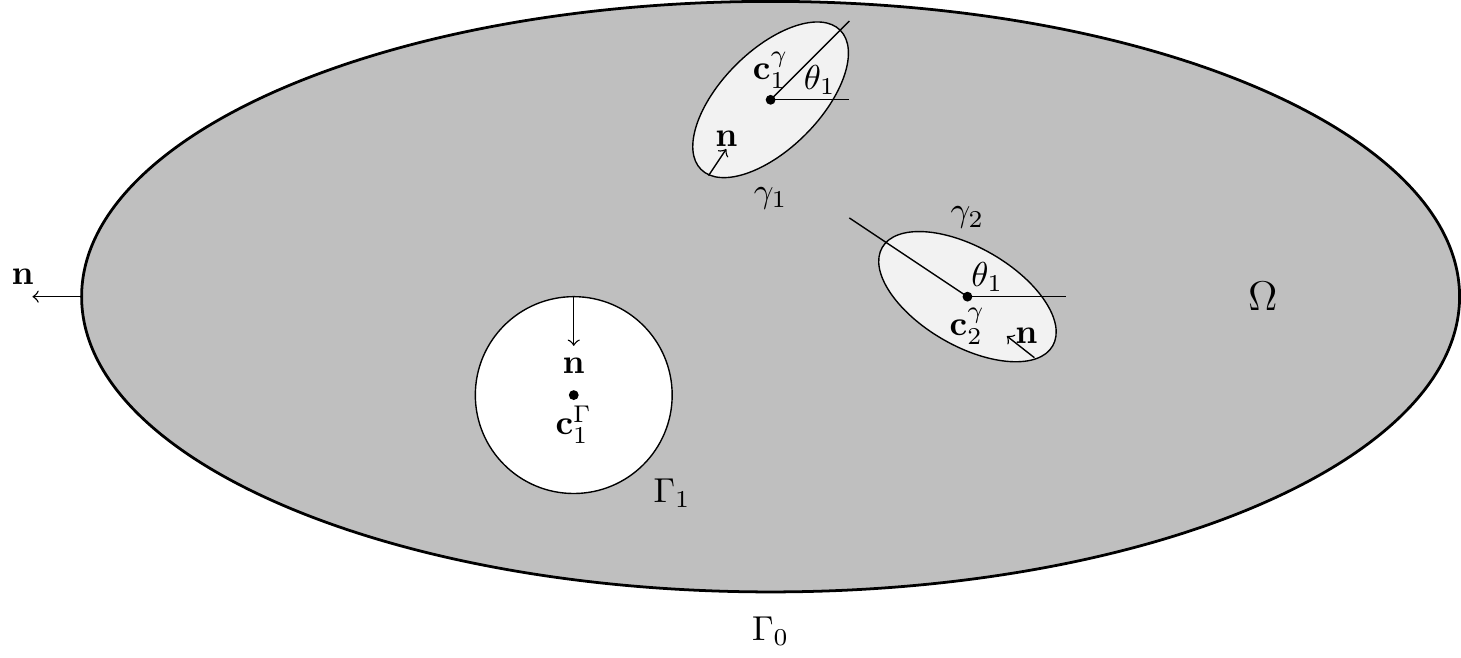}
\end{center}
\caption{\label{fig:geomSchematic}A sketch of a bounded fluid domain
$\Omega$.  $\gamma_1$ and $\gamma_2$ enclose rigid particles, while
$\Gamma_1$ is a solid wall.  If $\Omega$ is unbounded,  $\Gamma_0$  is
not present.  The vector $\nn$ is the unit normal vector pointing out of
the fluid domain.}
\end{figure}

We consider a collection of rigid particles suspended in a
two-dimensional bounded or unbounded domain, $\Omega$, with boundary
$\partial\Omega$. We let $\Gamma$ be the boundary of the fluid geometry,
$\Gamma_0$ is the outermost boundary if the domain is bounded, and
$\Gamma_i$, $1\leq i \leq M_w$ are the interior components of $\Gamma$.
The boundaries of rigid particles are $\gamma_j$, $1\leq j\leq M_p$, and
$\gamma = \cup_{j} \gamma_j$. Therefore, the fluid domain boundary is
$\partial\Omega =\Gamma \cup \gamma$, and we let $\nn$ be its outward
unit normal.  For each interior solid wall, we choose a single fixed
interior point $\cc^\Gamma_i$, and for each rigid particle, we require
an interior point $\cc^\gamma_i$ and a corresponding orientation angle
$\theta_i$.  A schematic of the geometry is in
Figure~\ref{fig:geomSchematic}.

\subsection{Governing Equations}\label{sec:governing}

We are interested in small particles and slow velocities which renders
the Reynolds number small $\Re \ll 1$, and the fluid is governed by the
incompressible Stokes equations.  A Dirichlet boundary condition $\UU$
is imposed on the solid walls $\Gamma$, a no-slip boundary condition is
imposed on the rigid bodies $\gamma$, and the rigid bodies are assumed
to be force- and torque-free.  On each solid wall, there is a net force
and torque, $\FF^\Gamma_i$ and $L^\Gamma_i$, respectively, that depend
on the boundary condition.  A similar force, $\FF^\gamma_j$, and torque
$\L^\gamma_j$ are defined for each rigid body $\gamma_j$, but these, for
the time being, are assumed to be 0.  Therefore, the governing equations
for $M_p$ particles suspended in a bounded $M_w$-connected domain is
\begin{equation}
  \label{eqn:modelEquations}
  \begin{split}
  \mu \Delta \uu = \grad p, &\hspace{20pt} \xx \in \Omega, \gap
    &&\mbox{\em conservation of momentum,}\\
  \grad \cdot \uu = 0, &\hspace{20pt} \xx \in \Omega, \gap
    &&\mbox{\em conservation of mass,} \\
  \uu = \UU, &\hspace{20pt} \xx \in \Gamma, \gap 
    &&\mbox{\em wall velocity,} \\
  \uu = \uu^\tau_j + \omega_j(\xx-\cc^\gamma_j)^\perp,&\hspace{20pt} 
    \xx \in \gamma, \gap &&\mbox{\em no-slip on the bodies,} \\
  \FF_j^\gamma = 0, &\hspace{20pt}j=1,\ldots,M_p, \gap 
    &&\mbox{\em force-free bodies,} \\
  L_j^\gamma = 0, &\hspace{20pt}j=1,\ldots,M_p, \gap 
    &&\mbox{\em torque-free bodies.}
  \end{split}
\end{equation}
Here, $\uu$ is the velocity, $p$ is the pressure, $\mu$ is the fluid
viscosity, $\uu^\tau_j$ and $\omega_j$ are the translational and
rotational velocities of rigid body $j$, respectively, and
$\FF_j^\gamma$ and $L_j^\gamma$ are the net force and torque of rigid
body $j$.  In the Stokes limit, the fluid viscosity sets
the time scale, and we assume it is one throughout the paper.  In the
case that the fluid domain is unbounded, the wall velocity equation is
replaced with the far-field condition
\begin{align*}
  \uu(\xx) = \uu_\infty(\xx), \quad |\xx| \rightarrow \infty.
\end{align*}
Upon solving for the translational and rotational velocities, the rigid
body centers and inclination angles $(\cc_j,\theta_j)$,
$j=1,\ldots,M_p$, satisfy
\begin{align}
\begin{split}
  \frac{d\cc_j}{dt} = \uu^\tau_j, \qquad 
  \frac{d\theta}{dt} = \omega_j.
\end{split}
\label{eqn:centersAngles}
\end{align}
Equations~\eqref{eqn:modelEquations} and~\eqref{eqn:centersAngles} govern
the dynamics of the rigid body suspensions, and their numerical solution
is a focus on this paper.

We have assumed that the rigid bodies are force- and torque-free.
However, when two rigid bodies are brought sufficiently close together,
numerical errors can cause the rigid bodies to unphysically intersect.
To avoid contact,  we will later relax the force- and torque-free
conditions to guarantee that numerical errors do not cause rigid bodies
to come into contact.  This idea is first described for
vesicle suspensions by Lu et al.~\cite{Lu2017} and we summarize the
method in Section~\ref{sec:repulsion}.

\subsection{Boundary Integral Equation Representation}
There exist many numerical methods for
solving~\eqref{eqn:modelEquations} such as level set
methods~\cite{Dou2007}, immersed boundary methods~\cite{Mittal2005},
dissipative particle dynamics~\cite{Pivkin2010}, smoothed particle
hydrodynamics~\cite{Polfer2016}, and lattice Boltzmann
methods~\cite{Ladd1994a, Ladd1994b}. However, because the fluid
equations are linear a boundary integral equation (BIE)
formulation~\cite{Pozrikidis1992} is possible.  BIEs have several
advantages including that only the interface has to be tracked, which
simplifies the representation of complex and moving geometries, and
high-order discretizations are straightforward.  We now reformulate
equation~\eqref{eqn:modelEquations} as a BIE.

We start by formulating the incompressible Stokes equations in the
absence of rigid bodies.  The {\em double-layer potential} is the
convolution of the stresslet with an arbitrary density
function~\cite{Ladyzhenskaya1963, Pozrikidis1992},
\begin{align}
  \label{eqn:dlp}
  \uu(\xx) = \DD[\eeta](\xx) = \frac{1}{\pi}\int_{\Gamma}
  \frac{\rr\cdot\nn}{\rho^2}\frac{\rr \otimes \rr}{\rho^2}
  \eeta(\yy)~\text{d}s_{\yy}, \quad \xx \in \Omega,
\end{align}
where $\rr = \xx - \yy$, $\rho=|\rr|$, and $\eeta$ is an unknown density
function defined on $\bd\Omega$.  The double-layer
potential~\eqref{eqn:dlp} satisfies the incompressible Stokes equations,
and the Dirichlet boundary condition $\UU$ is also satisfied if
$\eeta$ satisfies~\cite{Pozrikidis1992}
\begin{align}
  -\frac{1}{2} \eeta(\xx_0) + \DD[\eeta](\xx_0) = \UU(\xx_0), 
    \quad \xx_0 \in \Gamma.
  \label{eqn:secondKindBIE}
\end{align}
The double-layer potential cannot represent rigid body motions that
satisfy the incompressible Stokes equations.  Following Power and
Miranda~\cite{Power1987, Power1993}, this is resolved by introducing
point forces and torques due to each interior component of the geometry
$\Gamma_j$, and the strengths of these forces and torques are related to
the density function $\eeta$.  By introducing the velocity fields due to
a point force (Stokeslet) and a point torque (rotlet), both centered at
$\cc$,
\begin{align*}
  \mathbf{S}(\xx,\cc) = \frac{1}{4\pi}\left(-\log\rho\mathbf{I} + 
  \frac{\rr \otimes \rr}{\rho^2}\right), \quad \text{and} \quad
  \mathbf{R}(\xx,\cc) = \frac{\rr^\perp}{4\pi\rho^2},
\end{align*}
where $\rr = \xx - \cc$ and $\rho = |\rr|$.  Then, the second-kind
integral equation~\eqref{eqn:secondKindBIE} is replaced with the
completed second-kind BIE
\begin{equation}
  \label{eqn:completed_DLP}
  \begin{split}
  -\frac{1}{2}\eeta(\xx_0) + \DD[\eeta](\xx_0) + 
    \sum_{j=1}^{M_w} \left(\mathbf{S}(\xx,\cc^\Gamma_j)\FF^\Gamma_j + 
      \mathbf{R}(\xx,\cc^\Gamma_j)L^\Gamma_j\right) &= \UU(\xx_0),
      \quad &&\xx_0 \in \Gamma, \\
  \int_{\Gamma_j} \eeta~\text{d}s &= \FF^\Gamma_j, 
      &&j=1,\ldots,M_w, \\
  \int_{\Gamma_j} \eeta\cdot (\xx - \cc^\Gamma_j)^\perp~\text{d}s &=   
      L^\Gamma_j, &&j=1,\ldots,M_w.
  \end{split}
\end{equation}

We now introduce a suspension of rigid bodies $\gamma_j$,
$j=1,\ldots,M_p$.  The double-layer potential now includes contributions
from both the solid walls and rigid bodies.  Imposing the no-slip
boundary condition on the rigid bodies, a BIE formulation of the
suspension of rigid bodies governed by
equation~\eqref{eqn:modelEquations} is
\begin{subequations}
  \label{eqn:BIEformulation}
  \begin{align}
    \UU(\xx) &= -\frac{1}{2}\eeta(\xx) + \DD[\eeta](\xx) +
    \sum_{j=1}^{M_w} \left(\mathbf{S}(\xx,\cc^\Gamma_j)\FF^\Gamma_j + 
      \mathbf{R}(\xx,\cc^\Gamma_j)L^\Gamma_j\right)  \nonumber \\
&\hspace{96pt}+\sum_{j=1}^{M_p} \left(\mathbf{S}(\xx,\cc^\gamma_j)\FF^\gamma_j +
\mathbf{R}(\xx,\cc^\gamma_j)L^\gamma_j\right),
    \quad \xx \in \Gamma, \label{eqn:BIEformulation1} \\
  \uu^\tau_j + \omega_j(\xx - \cc_j^\gamma)^\perp &=
    -\frac{1}{2}\eeta(\xx) + \DD[\eeta](\xx) + 
    \sum_{j=1}^{M_w} \left(\mathbf{S}(\xx,\cc^\Gamma_j)\FF^\Gamma_j + 
      \mathbf{R}(\xx,\cc^\Gamma_j)L^\Gamma_j\right) \nonumber \\
&\hspace{96pt}+\sum_{j=1}^{M_p} \left(\mathbf{S}(\xx,\cc^\gamma_j)\FF^\gamma_j +
\mathbf{R}(\xx,\cc^\gamma_j)L^\gamma_j\right),
    \quad \xx \in \gamma, \label{eqn:BIEformulation2} \\
  \int_{\Gamma_j} \eeta~\text{d}s &= \FF^\Gamma_j, \quad
  \int_{\Gamma_j} \eeta\cdot (\xx - \cc^\Gamma_j)^\perp~\text{d}s =
  L^\Gamma_j, \quad j=1,\ldots,M_w, \label{eqn:BIEformulation3} \\
  \int_{\gamma_j} \eeta~\text{d}s &= \FF^\gamma_j, \quad
  \int_{\gamma_j} \eeta\cdot (\xx - \cc^\gamma_j)^\perp~\text{d}s =
  L^\gamma_j,\quad j=1,\ldots,M_p, \label{eqn:BIEformulation4} \\
  \FF^\gamma_j &= 0, \quad \L^\gamma_j = 0,\quad j=1,\ldots,M_p.
  \label{eqn:BIEformulation5}
\end{align}
\end{subequations}
Again, the methodology of Power and Miranda relates the strength of the
Stokeslets and rotlets of each rigid body to its density function.  The
BIE formulation~\eqref{eqn:BIEformulation} of the governing
equations~\eqref{eqn:modelEquations} consists of eight equations for
eight unknowns: the density function, net force, and net torque on the
solid walls and rigid bodies, and the translational and rotational
velocities.

While~\eqref{eqn:BIEformulation1} and~\eqref{eqn:BIEformulation2} are
both numerically desirable second-kind Fredholm integral equations
equations, equation~\eqref{eqn:BIEformulation1} has a rank one null
space because of the flux-free condition of the boundary data
$\UU$~\cite{Ladyzhenskaya1963}.  Following~\cite{Power1993}, this null
space is removed by adding the term 
\begin{align}
\label{eqn:N0_modification} 
  \mathcal{N}_0[\eeta](\xx) = \int_{\Gamma_0} 
    \nn(\xx)\otimes\nn(\yy)~\text{d}s(\yy)
\end{align}
to~\eqref{eqn:BIEformulation1}, but only for points $\xx \in \Gamma_0$.
Finally, if $\Omega$ is unbounded, equation~\eqref{eqn:BIEformulation}
has no null space, and the only modification is that
equation~\eqref{eqn:BIEformulation1} is removed and
equation~\eqref{eqn:BIEformulation2} has the background velocity
$\uu_\infty(\xx)$ is added to its right hand side.

\subsection{A Contact-Based Repulsion Force}
\label{sec:repulsion}
Exact solutions of the Stokes equations prohibit contact between force-
and torque-free bodies in finite time.  Therefore, any contact between
rigid bodies is caused by numerical errors.  The two main sources of
error are the quadrature error and time stepping error.  We address the
quadrature error using a combination of upsampling and interpolation,
and details of the method are described in~\cite{Quaife2014}.  Time
stepping methods have recently received a lot of attention, and some
recent works address stiffness~\cite{Quaife2014} and adaptive time
stepping~\cite{Kropinski1999, Quaife2015, Sorgentone2018}.  We describe
time stepping methods in Section~\ref{sec:temporal}.

We adopt the method of using artificial forces to avoid contact, and
there are many choices for the force. One possibility is a Morse or
Lennard-Jones potential that grows as a high order polynomial as two
bodies approach~\cite{Flormann2017, Liu2006}. This has been shown to
work for dense suspensions, but the resulting ODEs for the rigid body
dynamics become very stiff as the separation between bodies decreases.
Spring models~\cite{Tsubota2006, Zhao2013, Kabacogulu2017} have also
been used to generate artificial repulsion forces, but these models also
introduce stiffness. A further disadvantage of many contact algorithms
is that they do not guarantee that particles remain contact-free.

We apply a modification of the contact-aware method of Lu et
al.~\cite{Lu2017} that explicitly requires that particles remain
contact-free. We only summarize the method, but a in-depth description
is in~\cite{Lu2017}.  The method starts with the Stokes equations in
variational form,
\begin{align}
  \min \int_{\Omega} \nabla\uu:\nabla\uu~\text{d}\Omega,
  \quad\text{ such that }\quad \nabla\cdot\uu = 0, 
  \quad \xx \in \Omega.
  \label{eqn:stokesEL}
\end{align} 
After taking a time step, all pairs of bodies whose separation falls
below a minimum separation distance, which includes those that have
collided, are flagged.  Then, the space-time interference volume (STIV),
defined as the volume in space-time distance swept out by the trajectory
of two bodies, is computed.  Next, equation~\eqref{eqn:stokesEL} is
supplemented with an additional constraint that the STIV is positive (no
contact).  As detailed in~\cite{Lu2017, Harmon2011}, STIV offer a metric
to quantify collision volumes.  By using computing the STIV, rather than
simply the overlap at the new time step, time steps that result in rigid
bodies passing through each other are detected and resolved.  Details
for solving~\eqref{eqn:stokesEL} with the STIV inequality constraint are
given in Section~\ref{sec:contact}.

\subsection{Computing the Pressure and Stress}
To understand the rheological and statistical properties of suspensions
of rigid bodies, it is necessary to compute the pressure, stress, and
energy dissipation.  This expressions are standard~\cite{Power1993}, but
we summarize them for completeness.  A numerical method for Computing
the pressure, stress, and energy dissipation is described in
Section~\ref{sec:spatial}.

The pressure of the two-dimensional double-layer potential
is~cite{Power1993}
\begin{align*}
  p(\xx) = \frac{1}{\pi}\int_{\partial\Omega}\frac{1}{\rho^2}\left(
\mathbf{I} - 2\frac{\rr \otimes \rr}{\rho^2}\right)\nn\cdot\eeta~\text{d}s,
\qquad\xx\in\Omega.
\end{align*}
The pressure of the Stokeslet and rotlet contributions are easily
computed, and the pressure of the completed double-layer
potential~\eqref{eqn:completed_DLP} is
\begin{align}
  \label{eqn:pressure} 
  p(\xx) = \frac{1}{\pi}\int_{\partial\Omega}\frac{1}{\rho^2}\left(
    \mathbf{I} - 2\frac{\rr \otimes \rr}{\rho^2}\right)\nn\cdot\eeta~\text{d}s +
  \sum_{i=1}^{M_w}
  \frac{\FF_i^\Gamma\cdot(\xx-\cc_i^\Gamma)}{2\pi|\xx-\cc_i^{\Gamma}|^2}
  + \sum_{i=1}^{M_p}
    \frac{\FF_i^\gamma\cdot(\xx-\cc_i^\gamma)}{2\pi|\xx-\cc_i^\gamma|^2}.
\end{align}

Starting with the pressure~\eqref{eqn:pressure}, we compute the stress
tensor
\begin{align*} 
\ssigma = -p \mathbf{I} + \left(\nabla \uu + (\nabla\uu)^\intercal\right), \qquad
\xx\in\Omega.
\end{align*}
We decompose the stress into contributions from the double-layer
potential, Stokeslets, and rotlets
\begin{align}
  \label{eqn:stress}
  \ssigma = \ssigma_{\eeta} + \ssigma_{S} + \ssigma_{R},
\end{align}
where
\begin{align*}
  \ssigma_{\eeta}(\xx) &= \frac{1}{\pi}\int_{\partial\Omega}\left( 
    \frac{\nn\cdot\eeta}{\rho^2}\mathbf{I} 
    -8\frac{(\rr\cdot\nn)(\rr\cdot\eeta)(\rr\otimes\rr)}{\rho^6} 
    +\frac{(\rr\cdot\nn)(\rr\otimes\eeta + \eeta\otimes\rr)}{\rho^4} 
+ \frac{(\rr\cdot\eeta)(\rr\otimes\nn +
\nn\otimes\rr)}{\rho^4}\right)\text{d}s,\\
  \ssigma_S(\xx) &= -\sum_{i=1}^{M_w} 
    \frac{\FF_i^\Gamma\cdot(\xx - \cc_i^\Gamma)}{\pi|\xx-\cc_i^\Gamma|^2}
        (\xx - \cc_i^\Gamma)\otimes(\xx - \cc_i^\Gamma)  -
    \sum_{i=1}^{M_p}
    \frac{\FF_i^\gamma\cdot(\xx - \cc_i^\gamma)}{\pi|\xx-\cc_i^\gamma|^2}
        (\xx-\cc_i^\gamma )\otimes(\xx-\cc_i^\gamma),\\
\ssigma_R(\xx) &= -\sum_{i=1}^{M_w} \frac{L_i^\Gamma}{2\pi|\xx-\cc_i^\Gamma|^2}
((\xx-\cc_i^\Gamma)\otimes(\xx-\cc_i^\Gamma)^\perp +
    (\xx-\cc_i^\gamma)^\perp\otimes(\xx-\cc_i^\gamma))  \\
\qquad&-\sum_{i=1}^{M_p}
\frac{L_i^\gamma}{2\pi|\xx-\cc_i^\gamma|^2}
    ((\xx-\cc_i^\gamma )\otimes(\xx-\cc_i^\gamma)^\perp + 
    (\xx-\cc_i^\gamma)^\perp\otimes(\xx-\cc_i^\gamma)).
\end{align*}

These expressions hold inside $\Omega$. On the boundary there is a jump
in the pressure and stress, as reported in~\cite{Quaife2014}.  These
boundary terms are required for our near-singular integration scheme
described in Section~\ref{sec:spatial}.  Finally, using the divergence
theorem, we the volumed average stresses can be expressed in terms of
boundary integrals~\cite{Pozrikidis1992}. 
		
\section{Numerical Methods\label{s:method}} 
We are ultimately interested in studying statistical and rheological
properties of the fiber suspensions.  Therefore, we develop stable
numerical methods that solve the governing
equations~\eqref{eqn:modelEquations} for long time horizons.  To achieve
high-order accuracy, the spatial grid is resolved with spectral accuracy
(Section~\ref{sec:spatial}), and stability is achieved by using implicit
interactions in the time integrator and adaptive time stepping
(Section~\ref{sec:temporal}).  We are much more concerned with stability
than accuracy, so we use low-order but stable time stepping schemes.
Each time step requires solving a block diagonal preconditioned dense
linear system that is solved with GMRES and accelerated with the fast
multipole method~\cite{Greenbaum1992} (Section~\ref{sec:fast}).  We also
discuss a numerical methods for the collision algorithm
(Section~\ref{sec:contact}), and for computing the pressure and stress
(Section~\ref{sec:pressure_stress}).

\subsection{Spatial Discretization}\label{sec:spatial}
Let $\xx(\alpha)$, $\alpha \in [0,2\pi)$, be a parameterization of a
rigid body $\gamma_i$.  We will use a collocation method that requires
the discretization points $\xx_k$ at $k=1,\ldots,N_p$.  Spectral
accuracy is achieved by representing functions defined on $\gamma_i$ as
a Fourier series
\begin{align}
  f(\alpha) = f(\xx_i(\alpha)) = 
    \sum_{k = 0}^{N_p-1} \hat{f}_k e^{ik\alpha}.
\end{align}
The rigid walls are identically discretized at $N_w$ points and
functions defined on the rigid walls are also represented with a Fourier
series.  We use the FFT to compute the Fourier coefficients, and all
derivatives are computed with spectral accuracy using by differentiating
in Fourier space.

With collocation points defined on the rigid walls and solid bodies, we
now discretize the layer potentials.  We use a Nystr\"om method by
approximating the double-layer potentials with the trapezoid rule.
Equations~\eqref{eqn:BIEformulation} are discretized as
\begin{equation*}
  \begin{aligned}
  \UU(\xx_i) = -\frac{1}{2}\eeta(\xx_i) + 
  \sum_{k=1}^{N} K(\xx_i,\xx_k) \eeta(\xx_k) \Delta s_k
    &+ \sum_{j=1}^{M_p} \left(\mathbf{S}(\xx_i,\cc^\gamma_j)\FF_j +
    \mathbf{R}(\xx_i,\dd_j)L_j\right)  \\
    &+ \sum_{j=1}^{M_w} \left(\mathbf{S}(\xx_i,\cc^\Gamma_j)\FF_j +
    \mathbf{R}(\xx_i,\cc^\Gamma_j)\L_j\right), \quad \xx_i \in \Gamma,
  \end{aligned}
\end{equation*}
for $\xx_i \in \Gamma$, and
\begin{equation*}
  \begin{aligned}
\uu^\tau_j + \omega_j(\xx_i - \cc^\gamma_j)^\perp =
-\frac{1}{2}\eeta(\xx_i) +
\sum_{k=1}^{N} K(\xx_i,\xx_k) \eeta(\xx_k) \Delta s_k
    &+ \sum_{j=1}^{M_p} \left(\mathbf{S}(\xx_i,\cc^\gamma_j)\FF_j +
    \mathbf{R}(\xx_i,\dd_j)L_j\right) \\
    &+ \sum_{j=1}^{M_w} \left(\mathbf{S}(\xx_i,\cc^\Gamma_j)\FF_j +
    \mathbf{R}(\xx_i,\cc^\Gamma_j)\L_j\right),
  \end{aligned}
\end{equation*}
for $\xx \in \gamma_j$, where $N = M_w N_w + M_p N_p$ is the total
number of discretization points and
\begin{align*}
  K(\xx,\yy) = \frac{1}{\pi} \frac{\rr \cdot \nn}{\rho^2} 
               \frac{\rr \otimes \rr}{\rho^2}
\end{align*}
is the kernel of the double-layer potential.  The diagonal entries of
$K$ are replaced with the limiting value
\begin{align*}
  \lim_{\substack{\yy \rightarrow \xx \\ \yy \in \bd\Omega}} 
    K(\xx,\yy) = \frac{\kappa}{2\pi}\tt\otimes\tt,
    \quad \xx \in \bd\Omega,
\end{align*}
where $\kappa(\xx)$ is the curvature and $\tt(\xx)$ is the tangent
vector of $\bd\Omega$ at $\xx$.

Since the kernel $K(\xx,\yy)$ is smooth, the trapezoid rule guarantees
spectral accuracy~\cite{Trefethan2014}.  However, at a fixed resolution
$N$, the error grows when a target point and a source point on different
bodies are sufficiently close.  The error is caused by a nearly-singular
integrand with large derivatives.  To resolve this issue, an algorithm
for near-singular integration method must be employed.  We use the
interpolation scheme outlined in~\cite{Quaife2014} that is based off of
the algorithm first outlined in~\cite{Ying2006}.
  
After discretizing~\eqref{eqn:BIEformulation} with a Nystr\"om method,
the result is a dense linear system for the density function, rotlets
and Stokeslets, and the translational and rotational velocities of each
body.  We chose a double-layer potential formulation so that the linear
system can be solved with a mesh-independent number of GMRES
iterations~\cite{Campbell1996}.  Therefore, the algorithmic cost is
dominated by the cost of a matrix-vector multiplication, and the number
of required geometry-dependent GMRES iterations.  Algorithms for
controlling these costs are discussed in Section~\ref{sec:fast}.

\subsection{Contact Resolution}
\label{sec:contact}
The contact resolution method starts by advancing the force- and
torque-free bodies from time $t$ to $t + \Delta t$ with
equation~\eqref{eqn:centersAngles}.  Then, equation~\eqref{eqn:stokesEL}
needs to be solved with the inequality constraint that the STIV is
greater than zero.  The definition of the STIV includes a parameter that
guarantees that not only do bodies intersect, but they maintain minimum
separation distance. Following~\cite{Lu2017}, we use a Lagrange
multiplier, $\lambda$, to satisfy the non-negative STIV inequality
constraint (a negative STIV indicates contact).  The resulting equations
are  the incompressible Stokes equations with body forces and torques
that depend on the gradient of the STIV, and the inequality constraint
that the STIV is less than or equal to zero.  Therefore, the governing
equations are similar to~\eqref{eqn:BIEformulation}, except that the
force- and torque-free conditions have been changed.

The resulting equation is a non-linear complementary problem (NCP) To
solve the NCP, the problem is linearized as a sequence of linear
complementary problems (LCPs) that are solved until the STIV is less
than or equal to zero.  For small minimum separation distances, it is
possible that the LCP converges very slowly, or not at all, to the
solution of the NCP.  In Section~\ref{sec:temporal}, we describe two
methods to avoid this slow convergence.

If only two bodies are in contact, the symmetry of the STIV ensures that
the net force is zero.   However, if more than two bodies are in
contact, then the forces returned by the STIV are not guaranteed to sum
to zero.  We improve the validity of the method by always requiring that
the total force added to the system is zero, meaning
\begin{align*}
  \sum_{j=1}^{M_p} \FF^\gamma_j = \mathbf{0}. 
\end{align*}
To accomplish this, we group all bodies that are in contact into
distinct clusters.  For each cluster, the body that is in contact with
the most other bodies (to break ties, the body with the largest
repulsion force in magnitude) is given a net force that balances all
other bodies in the cluster. Whenever a force on a body is scaled the
torque on that body is scaled by the same amount.

\subsection{Time Stepping Methods}
\label{sec:temporal}
Since we only consider rigid body suspensions, we only track the centers
$\cc_i^\gamma$ and orientations $\theta_i$ of each rigid body
$\gamma_i$. Given a suspension of rigid
bodies,~\eqref{eqn:BIEformulation} is solved for the translational and
rotational velocities of each body. Then, the position and angle of each
body are updated according to the ODEs,
\begin{align*}
  \frac{\text{d}}{\text{d}t}\cc_k &= \uu^\tau_k,  
    \quad \cc_k(0) = \cc_k^0, \\
  \frac{\text{d}}{\text{d}t}\theta_k &= \omega_k,
    \quad \theta_k(0) = \theta_k^0.
\end{align*}
The ODEs are advanced in time using the first-order explicit Euler
method
\begin{align*}
  \cc_k^{N+1} &= \cc_k^N + \Delta t \left(\uu_k^\tau\right)^N, \\
  \theta_k^{N+1} &= \theta_k^N + \Delta t \left(\omega_k\right)^N.
\end{align*}
A low cost first-order time stepper is justified since we are interested
in statistical, rheological, and bulk quantities rather than the
individual trajectories of the rigid bodies.  If high-order accuracy is
desired, a Runge-Kutta or deferred correction method~\cite{Quaife2015,
qua-bir2016} can be applied.

Since the dynamics of the suspension can develop complex features, we
use an heuristic adaptive time stepping method.  When the LCP solver
requires many iterations to solve the NCP, this indicates that the time
step size should be reduced.  Conversely, if the LCP solver converges
quickly to the solution of the NCP, then a larger time step size can be
taken.  Therefore, if the each LCP iteration is not converging
to a contact-free configuration sufficiently fast, we abort the time
step and restart with the time step size $\Delta t/2$. Conversely, if
the LCP iteration converges quickly, then, at the next time step, we
increase $\Delta t$ to $1.5\Delta t$.

The stability of the method depends on the discretization
of~\eqref{eqn:BIEformulation} and the minimum separation distance.  If
the minimum separation distance is small and the interaction between two
nearly-touching bodies is discretized explicitly, then convergence of
the LCP solver to the solution of the NCP requires a large number of
iterations, and the maximum stable time step size is small.  In
contrast, if the interaction between the bodies is discretized
implicitly, then the stiffness caused by the nearly-touching bodies is
resolved and a larger stable time step can be taken.

The leading source of stiffness for a particular rigid body is the
velocity field induced by the body itself. This can be controlled with a {\em locally implicit} discretization
\begin{align}
  \DD[\eeta](\xx) \approx \DD^{n}[\eeta^{N+1}](\xx) +
  \DD^{n}[\eeta^N](\xx),
  \label{eqn:BlockImplicit}
\end{align}
where the superscript of $\DD$ implies that the geometry used in the
discretization of the layer potential is at time step $n$.  This time
stepping method is focus of the work of Lu et al.~\cite{Lu2017}.
There, to avoid stiffness caused by nearly-touching bodies, a
sufficiently large minimum separation distance is used.
 
Instead of avoiding nearly-touching bodies, we use the {\em globally
implicit} discretization
\begin{align}
  \DD[\eeta](\xx) \approx \DD^{N}[\eeta_{N+1}](\xx) +
  \DD^{N}[\eeta^{N+1}](\xx).
  \label{eqn:SemiImplicit}
\end{align}
By discretizing the interactions between different bodies implicitly,
the stiffness caused by nearly-touching bodies is reduced, so large time
steps with small minimum separation distances are possible.

\subsection{Fast Summation and Preconditioning}
\label{sec:fast}
A discretization of the locally implicit time
stepper~\eqref{eqn:BlockImplicit} results in a block-diagonal linear
system, where each block is a dense $N_p \times N_p$ or $N_w \times N_w$
matrix.  In contrast, a discretization of the globally implicit time
stepper~\eqref{eqn:SemiImplicit} results in a dense $N \times N$ linear
system.  Therefore, considering only the rigid body contributions, the
cost of a single matrix-vector multiplication is $\mathcal{O}(M_p
N_p^2)$ for~\eqref{eqn:BlockImplicit} and $\mathcal{O}(M_p^2 N_p^2$
for~\eqref{eqn:SemiImplicit}.  We reduce the cost of matrix-vector
multiplication by using the fast multipole method.  This reduces the
cost of a single matrix-vector multiplication for both the locally
implicit and globally implicit time stepping methods to $\mathcal{O}(N)$
operations.

Reducing the cost of matrix-vector multiplication greatly reduces the
computational effort.  However, for dense suspensions, the number of
GMRES iterations can also be large.  Therefore, we apply a
block-diagonal preconditioner where each block is precomputed and
factorized.  The identical preconditioner has been used for
two-dimensional vesicle suspensions~\cite{Quaife2014}.  For rigid
bodies, the preconditioner is further accelerated by factorizing all
double-layer potentials at the initial condition, and then multiplying
with rotation matrices corresponding to the rotation velocity $\omega$.

\subsection{Computing the Pressure and Stress}
\label{sec:pressure_stress}
Computing the pressure~\eqref{eqn:pressure} and
stress~\eqref{eqn:stress} are more challenging than evaluating the
velocity double-layer potential.  The challenge stems from a singularity
that scales as $\mathcal{O}(|\xx - \yy|^{-2})$ as a target point
approaches a source point.  The pressure and stress are computed using a
combination of singularity subtraction and odd-even
integration~\cite{sid-isr1988}.  The result is a spectrally accurate
method for computing the pressure and stress.

Near-singular integration still required, and our interpolation-based
method requires limiting values of the layer
potentials~\eqref{eqn:pressure} and~\eqref{eqn:stress} as $\xx
\rightarrow \bd\Omega$.  The limiting values and tests of the
near-singular integration scheme are reported in~\cite{Quaife2014}.

\section{Results\label{s:results}} 
We use our new time stepping method to simulate bounded and unbounded
suspensions of two-dimensional rigid bodies in a viscous fluid.  The
main parameters are the minimum separation distance $\delta$, the number
of discretization points of each body, $N_p$, and each solid wall,
$N_w$, and the initial time step size $\Delta t$. We perform convergence
studies and investigate the effect of the STIV algorithm on the
reversibility of the flow.  To further demonstrate the consequence of
STIV, we include plots of streamlines that cross whenever the collision
detection algorithm is applied.  The particular experiments we perform
are now summarized.
\begin{itemize}
  \item {\bf Shear Flow}: We consider the standard problem of two
  identical rigid circles in the shear flow $\uu = (y,0)$ with the left
  body slightly elevated from the right body.  We report similar results
  to those presented in~\cite{Lu2017}, but we are able to take smaller
  initial displacements and minimum separation distances.  The contact
  algorithm breaks the reversibility of the flow, and this effect is
  illustrated and quantified.

  \item {\bf Taylor-Green Flow}: We simulate a concentrated suspension
  of 48 rigid ellipses in an unbounded Taylor-Green flow. At the
  prescribed separation distance, our new time stepping method is able
  to stably reach the time horizon, while the locally semi-implicit time
  integrator proposed by Lu et al.~\cite{Lu2017} results in the STIV
  algorithm stalling, even with $\Delta t = 10^{-8}$.

  \item {\bf Porous Monolayer Injection}: We consider a suspension of
  confined rigid circular bodies motivated by an experiment by MacMinn
  et al.~\cite{MacMinn2015}.  The geometry is an annulus with an inflow
  at the inner boundary and an outflow at the outer boundary.  We again
  examine the effect STIV on the reversibility of the flow, and compute
  the shear strain rate and make qualitative comparisons to results for
  deformable bodies~\cite{MacMinn2015}.

  \item {\bf Taylor-Couette Flow}: With the ability to do high area
  fraction suspensions without imposing a large non-physical minimum
  separation distance, we simulate rigid bodies of varying aspect ratios
  inside a Taylor-Couette device.  We examine the effect of the rigid
  body shape and area fraction on the effective viscosity and the
  alignment angles.
\end{itemize}

\subsection{Shear Flow}
\label{sec:shear}
We consider two rigid circular bodies in the shear flow ${\uu}(\xx) =
(y,0)$.  One body is centered at the origin, while the other body is
placed to the left and slightly elevated of the origin.  With this
initial condition, the particles come together, interact, and then
separate.  Both bodies are discretized with $N=32$ points and the arc
length spacing $h = 2\pi/32 \approx 0.196$.  This experiment was also
performed by Lu et al.~\cite{Lu2017}, and we compare the two time
stepping methods.

\begin{figure}[!h]
  \begin{center}
    \includegraphics{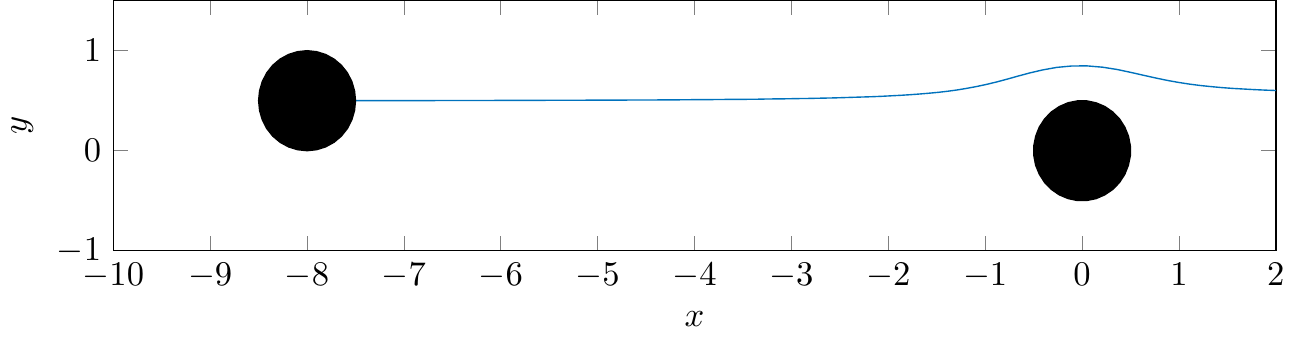}
    \begin{tabular}{c c}
      \includegraphics{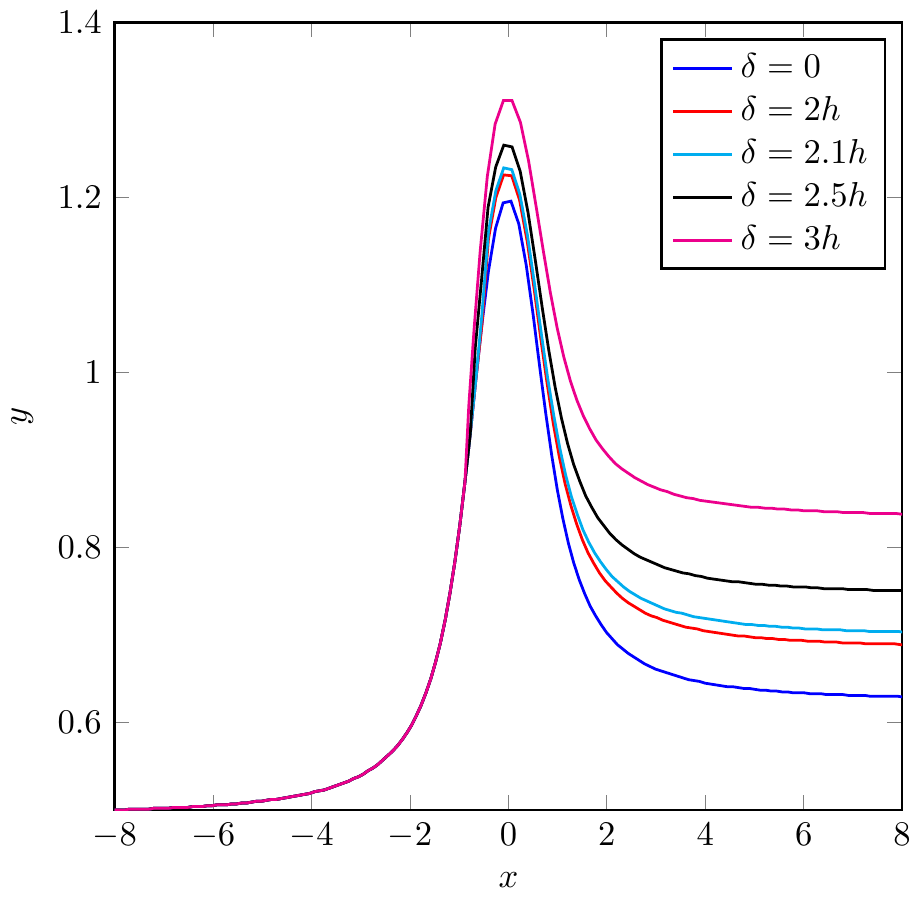} &
      \includegraphics{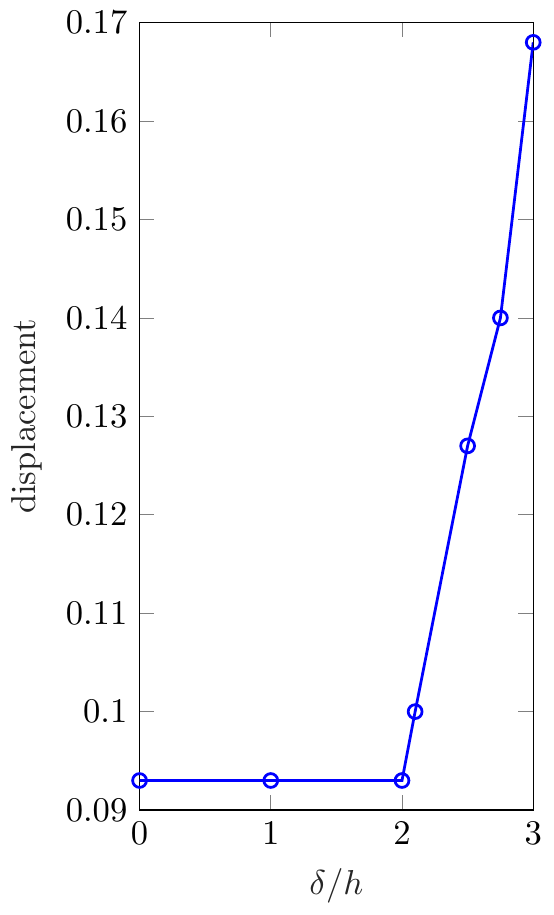}
    \end{tabular}
  \end{center}
\caption{\label{fig:shear_experiment} Shear experiment. Top: The initial
setup and trajectory of the left body.   Bottom left: The left body's
trajectory for varying minimum separation distances.  Notice how
the trajectories are identical until shortly before $x=0$ when the
contact algorithm is first applied.  Bottom right: The final vertical
displacement of the left particle for varying minimum separation
distances.}
\end{figure}

We start by considering the time step size $\Delta t = 0.4$ and
minimum separation distance $\delta = 0$ (no contact algorithm). 
Our new globally implicit method successfully reaches the time
horizon without requiring a repulsion force.  However, with the same
$\Delta t$, the local explicit time stepping results in a collision
between the bodies, so the collision algorithm is required to reach the
time horizon.  Alternatively, the time step size can be reduced, but, as
we will see, for sufficiently dense suspensions, even an excessively
small time step size results in collisions.  Next, in
Figure~\ref{fig:shear_experiment}, we investigate the effect of the
minimum separation distance on the position of the rigid bodies.  The
top plot shows the trajectory of the left body as it approaches,
interacts, and finally separates from the body centered at the origin.
In this simulation, we use our new globally implicit time integrator,
but the STIV contact algorithm is not applied.  The bottom left plot
shows the trajectory of the particle when the contact algorithm is
applied with varying levels of separation.  Notice that the trajectories
are identical until near $x=0$ when the particle separation first falls
below the minimum separation distance.  Finally, in the bottom right
plot, the final vertical displacement body initially on the left is
plotted.  These results are computed for the locally implicit time
stepping method~\cite{Lu2017}, and the general trend of the trajectories
are similar.

We next investigate the effect of the collision algorithm on the time
reversibility of the flow. We reverse the shear direction at $t=10$ and
measure the error between the body's center at $t=0$ and $t=20$.  We
expect an error that is the sum of a first-order error caused by time
stepping, and a fixed constant caused by the minimum separation
distance.  The results for various values of $\delta$ are reported in
Table~\ref{tab:reverse}. When the contact algorithm is not applied
when $\delta=0$ and $\delta=h$, we observe the expected
first-order convergence.  When $\delta \geq 2h$, the bodies are
deflected onto contact-free streamlines when their proximity reaches the
minimum separation distance.  After the flow is reversed, the bodies
again pass one another, but they are now on contact-free streamlines, so
the initial deflection is not reversed.  For these larger values of
$\delta$, we see in Table~\ref{tab:reverse} that the error eventually
plateaus as $\Delta t$ is decreased.

\begin{table}[!h]
\begin{center}
\begin{tabular}{c| c c c c c}
$ $ & & & $\Delta t$ & &\\
$\delta$ & $4\e{-2}$ &$ 2\e{-2}$ & $1\e{-2}$ & $5\e{-3}$ & $2.5\e{-3}$\\
\hline
0 & $1.35\e{-1}$ & $7.32\e{-2}$ & $3.74\e{-2}$ & $2.00\e{-2}$ & $1.01\e{-2}$\\
$h$ & $1.35\e{-1}$ & $7.32\e{-2}$ & $3.74\e{-2}$ & $2.00\e{-2}$ & $1.01\e{-2}$\\
$2h$ & $1.88\e{-1}$ & $1.41\e{-1}$ & $1.17\e{-1}$ & $1.08\e{-1}$ &
$1.02\e{-1}$\\
$2.25h$ & $2.55\e{-1}$ & $2.08\e{-1}$ & $1.87\e{-1}$ & $1.78\e{-1}$ &
$1.73\e{-1}$\\
$2.50h$ & $3.05\e{-1}$ & $2.69\e{-1}$ & $2.52\e{-1}$ & $2.45\e{-1}$ &
$2.40\e{-1}$\\
$2.75h$ & $3.64\e{-1}$ & $3.31\e{-1}$ & $3.13\e{-1}$ & $3.07\e{-1}$ &
$3.03\e{-1}$\\
$3.00h$ & $4.12\e{-1}$ & $3.88\e{-1}$ & $3.72\e{-1}$ & $3.67\e{-1}$ &
$3.63\e{-1}$
\end{tabular}
\end{center}
\caption{A study of time reversibility of the shear flow example. At
$t=10$, the flow direction is reversed and we calculate the relative
error in the initial and final positions. When the collision constraint
is active and force is needed to keep the bodies apart the error in
the reversibility is dominated by the contact
algorithm.}\label{tab:reverse}
\end{table}

\begin{figure}[!h]
\begin{center}
\includegraphics{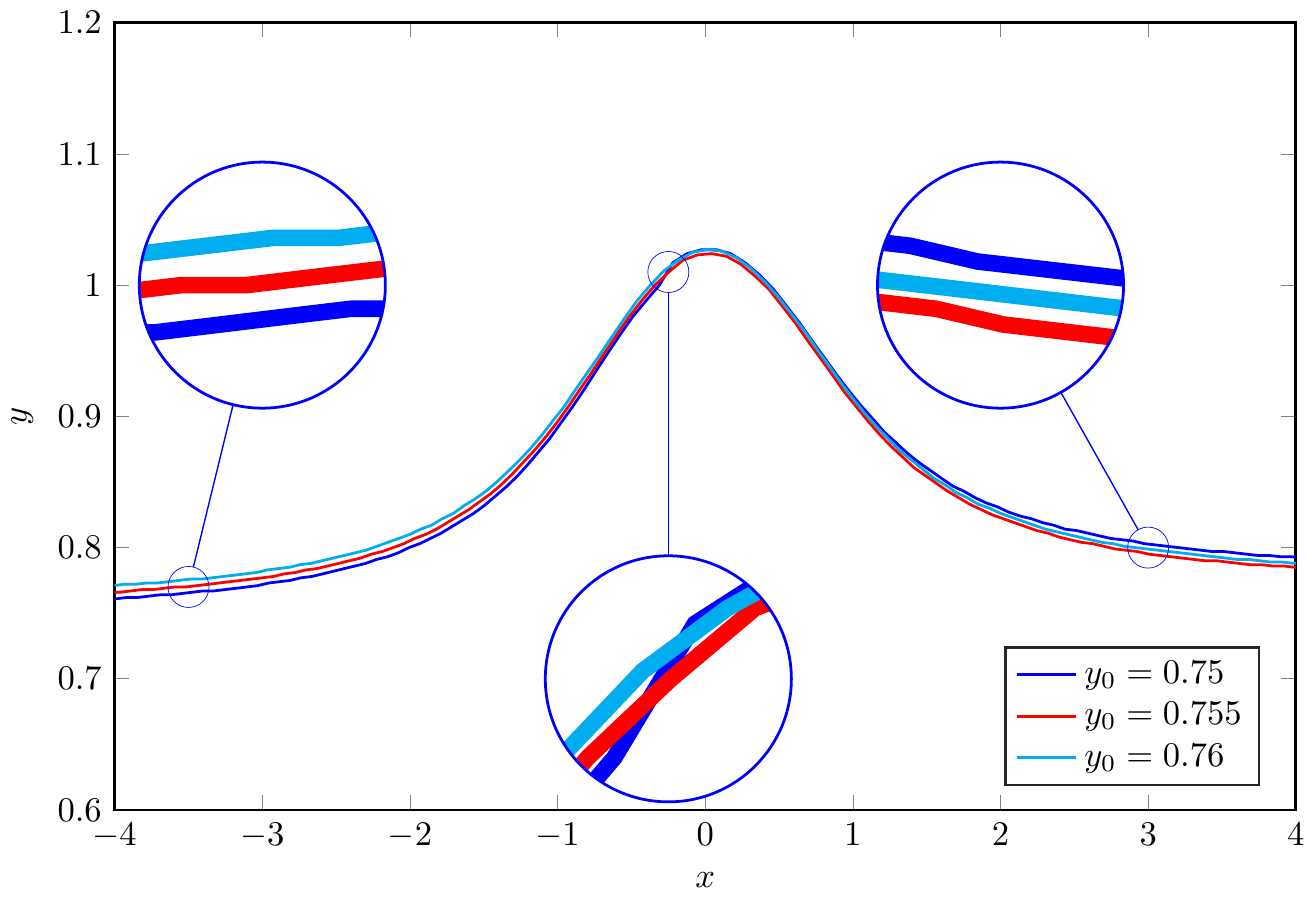}
\end{center}
\caption{\label{fig:shear_cross} The contact algorithm causes
streamlines to cross. Keeping the minimum separation fixed at
$\delta=3h$, we vary the starting $y$ location of the left body. The
teal and red streamlines do not require a repulsion force to enforce the
minimum separation between the bodies, but the blue streamline does.
Once the contact algorithm is applied, the blue streamline crosses the
other streamlines (middle inset). This crossing of the streamlines
breaks the reversibility of the simulation.}
\end{figure}

The break in reversibility is further demonstrated by examining
individual streamlines.  In Figure~\ref{fig:shear_cross}, we compute the
streamline of the left body for three different initial placements.  We
set $\delta=3h$ for all the streamlines.  With this threshold, only the
bottom-most streamline falls below $\delta$.  Therefore, as the bodies approach, the streamlines behave as expected---they do
not cross.  However, when the contact algorithm is applied to the blue
streamline, the streamlines cross.

\subsection{Taylor-Green Flow}

For planar flows, we can separate suspensions into dilute and
concentrated regimes by comparing the number of bodies per unit area,
$\nu$, to the average body length $\ell$. If $\nu < 1/\ell^2$, then we
are in the dilute regime, otherwise we are in the concentrated regime
(in 2D planar suspensions, unlike 3D suspensions, there is no
semi-dilute regimes).  We consider the suspension of 75 rigid bodies in
the Taylor-Green flow $\mathbf{u}^\infty = (\cos(x)\sin(y),
-\sin(x)\cos(y))$.  The number of bodies per unit area is $\nu \approx
3.1$ which is greater than $1/\ell^2=1.1$.  Therefore, this suspension
is well within the concentrated regime. 

\begin{figure}[!h]
  \begin{center}
    \begin{tabular}{c c }
      \includegraphics[width=6cm]{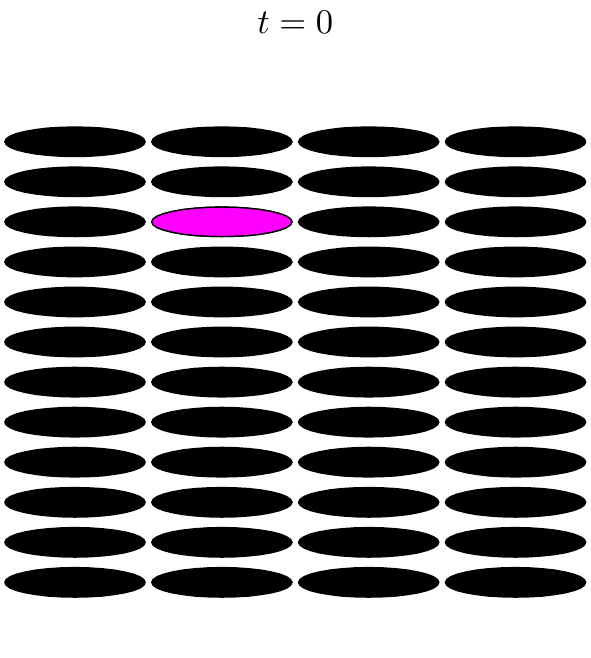} &
      \includegraphics[width=6cm]{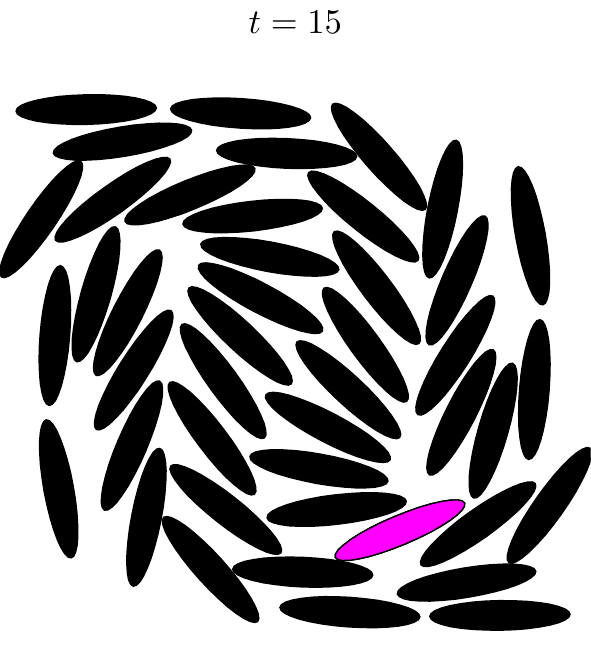}\\
      \includegraphics[width=6cm]{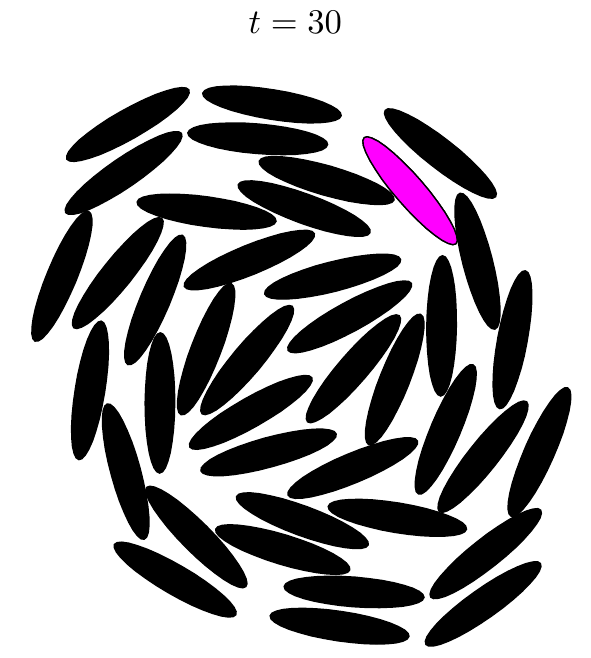} &
      \includegraphics[width=6cm]{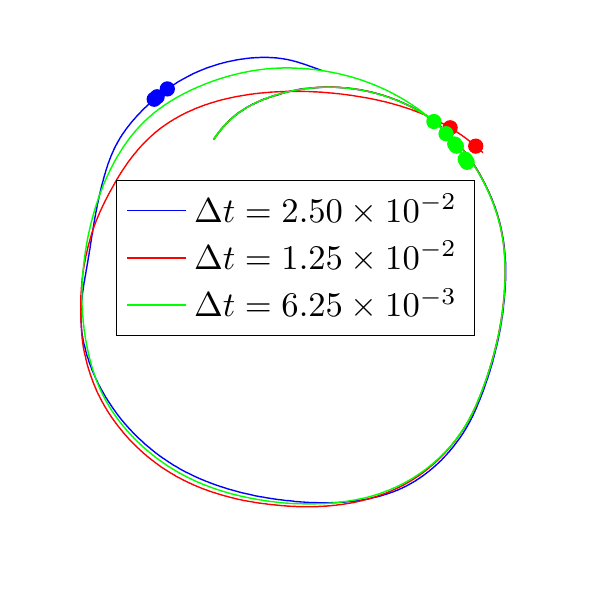}
    \end{tabular}
  \end{center}
  \caption{\label{fig:taylor_green} Snapshots of a dense suspension in
  an unbounded Taylor-Green flow.  The number of bodies per unit area
  $\nu$ is approximately 3.1. This is greater than $1/\ell^2 = 1.1$,
  which puts the simulation well within the concentrated regime.  Bodies
  are discretized with 32 points and the minimum separation is
  $\delta=0.05h$.  The bottom right plot shows the trajectory of the
  center of the colored body for different step sizes. Each line in that
  plot is marked where a repulsion force is used to enforce the minimum
  separation. }
\end{figure}

We discretize the bodies with $N=32$ points and select the minimum
separation distance $\delta=0.05h$. Snapshots of the simulation are
shown in Figure \ref{fig:taylor_green}.  In this concentrated
suspension, the bodies come into contact much more frequently.  If
the interactions between these nearly touching bodies are treated
explicitly, this leads to stiffness.  Our time stepper controls this
stiffness by treating these interactions implicitly, and the simulation
successfully reaches the time horizon.  We performed the same
simulation, but with the locally implicit time stepping
method~\cite{Lu2017}.  Because of the near-contact, smaller time step
sizes must be taken.  We took time step sizes as small as $10^{-8}$, and
the method was not able to successfully reach the time horizon.  This
exact behavior has also been observed for vesicle
suspensions~\cite{Quaife2014}.  In the bottom right plot of
Figure~\ref{fig:taylor_green}, we show the trajectory of one body for
different time step sizes.  The dots denote locations where the contact
algorithm is applied.  For this very complex flow, the trajectories are
in good agreement with different time step sizes.

\FloatBarrier
\subsection{Fluid Driven Deformation}
\begin{figure}[h!]
\begin{center}
\includegraphics{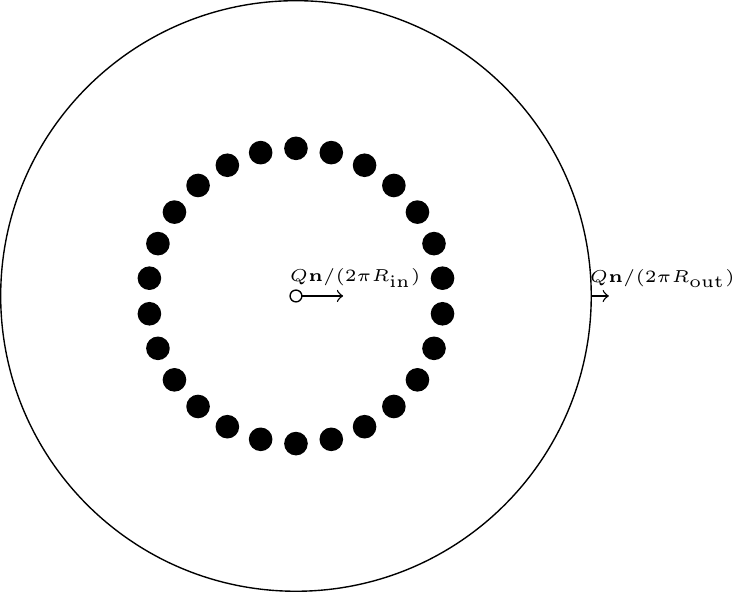}
\end{center}
\caption{\label{fig:radial} The geometry used in our numerical
experiment that is motivated by the experimental setup of MacMinn et
al.~\cite{MacMinn2015}.  The fixed solid bodies are shaded in black.}
\end{figure}

A recent experiment considers a dense monolayer packing of soft
deformable bodies~\cite{MacMinn2015}.  Motivated by this experiment, we
perform numerical simulations of rigid bodies in a similar device.  We
pack rigid bodies in a Couette device, but with a very small inner
boundary.  The boundary conditions are an inflow or outflow of rate $Q$
at the inner boundary with an outflow or inflow at the outer cylinder.
This boundary condition corresponds to injection and suction of fluid
from the center of the experimental microfluidic device.  In the
experimental setting, the soft bodies are able to reach the outer
boundary, and the resulting boundary condition would not be uniform at
the outer wall.  So that we can apply the much simpler uniform inflow or
outflow at the outer boundary, we force the rigid bodies to remain
well-separated from the outer wall.  We accomplish this by placing a
ring of {\em fixed} rigid bodies halfway between the inner and outer
cylinders (Figure~\ref{fig:radial}).  The spacing between these fixed
bodies is sufficiently small that the mobile bodies are not able to
pass.  Since the outer boundary is well-separated from the fixed bodies,
the outer boundary condition is justifiably approximated with a uniform
flow.

\begin{figure}[h!]
  \begin{tabular}{c c c}
    \includegraphics{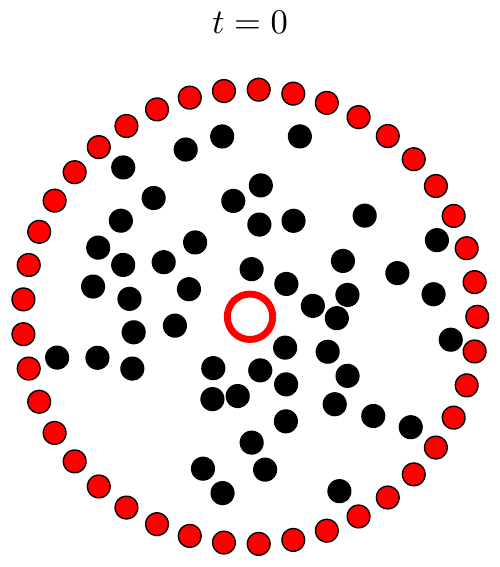}&
    \includegraphics{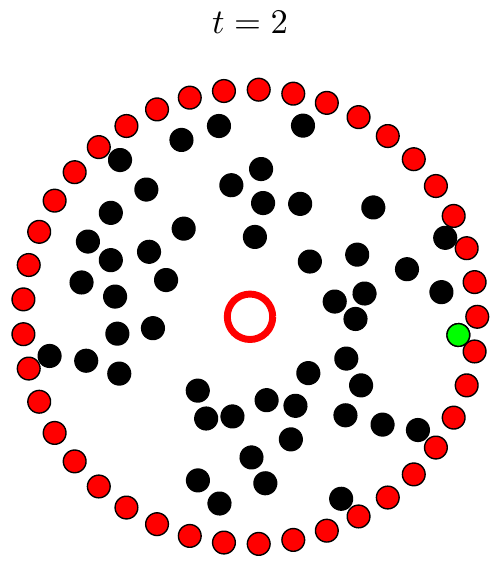}&
    \includegraphics{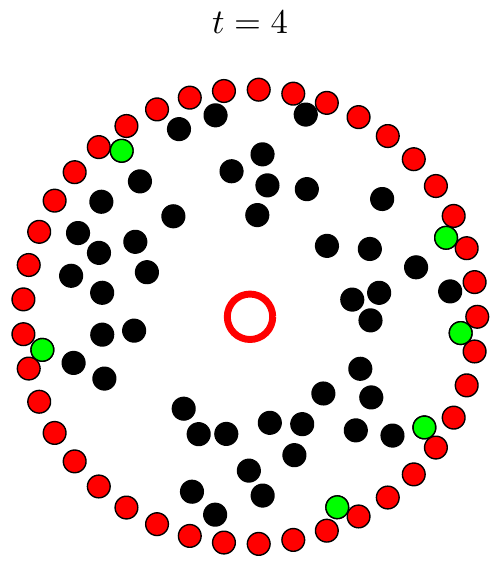}\\
    \includegraphics{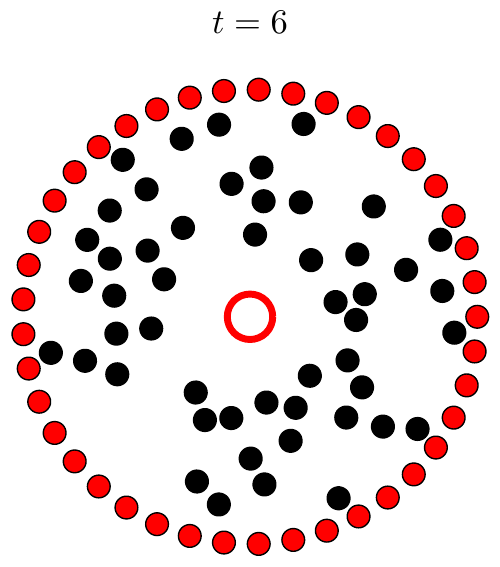}&
    \includegraphics{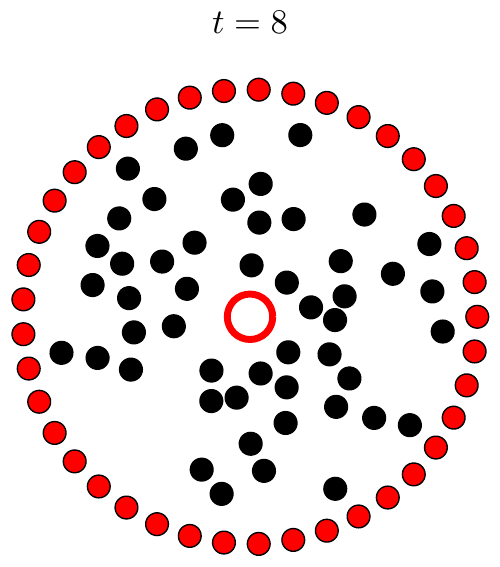}&
    \includegraphics{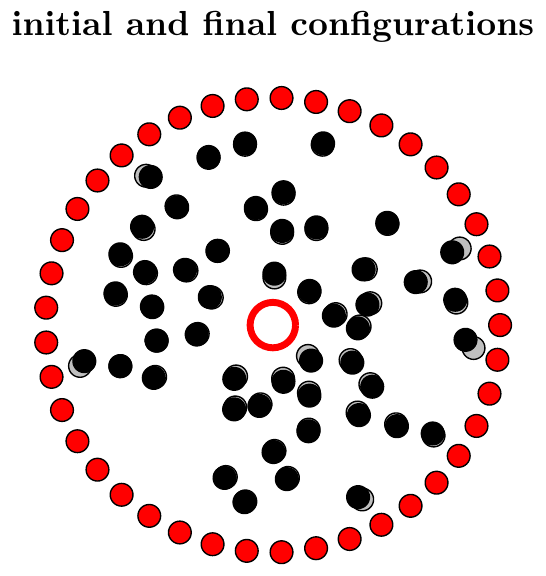}\\
  \end{tabular}
  \caption{\label{fig:macminn} Snapshots of a rigid body suspension
  motivated by an experiment for deformable bodies~\cite{MacMinn2015}.
  Fluid is injected at a constant rate starting at $t=0$. At $t=4$ the
  flow direction is reversed. Fixed bodies are colored in red, while
  bodies subject to a repulsion force are colored in green. The initial
  configuration has been superimposed on the final configuration at
  $t=8$ to show the effect of the repulsion forces on reversibility.}
\end{figure}

We start by examining the effect of the contact algorithm on the
reversibility of the flow.  We again reverse the flow at time $T$ and
run the simulation until time $2T$. The rigid bodies are in contact for
much longer than the shear example in Section~\ref{sec:shear}, so
maintaining reversibility is much more challenging.
Figure~\ref{fig:macminn} shows several snapshots of the simulation, and
the bottom right plot superimposes the initial and final configurations.
We observe only a minor violation of reversibility, and it is largest
for bodies that were initially near the fixed bodies---the contact
algorithm is applied to these bodies most frequently.

\begin{figure}[h!]
  \begin{tabular}{c c c}
    \includegraphics{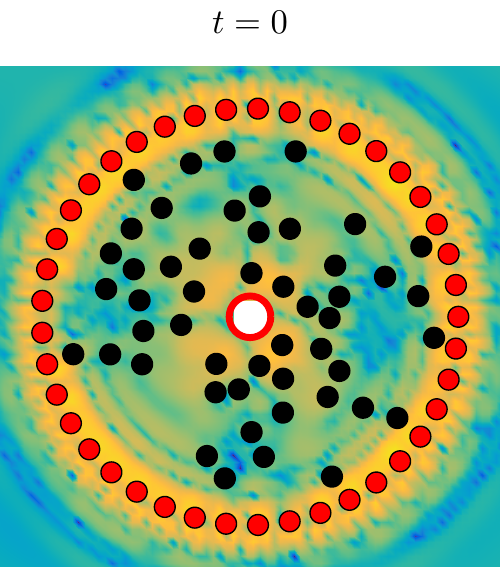}&
    \includegraphics{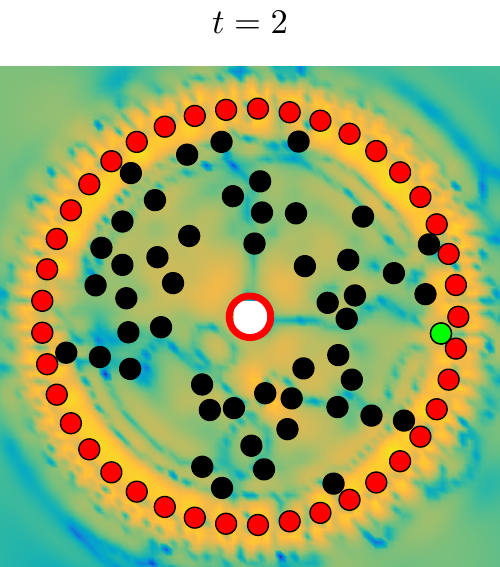}&
    \includegraphics{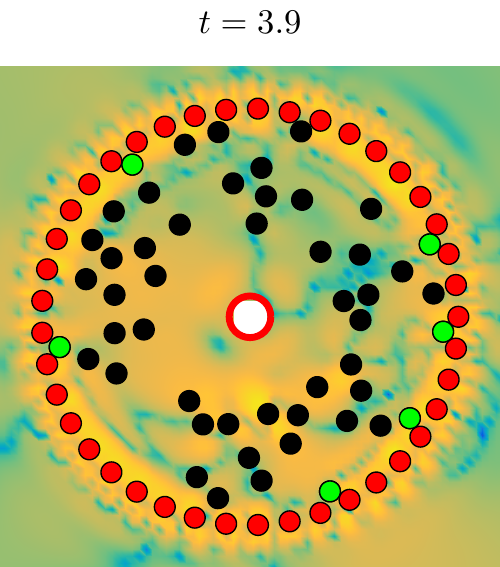}\\
    \includegraphics{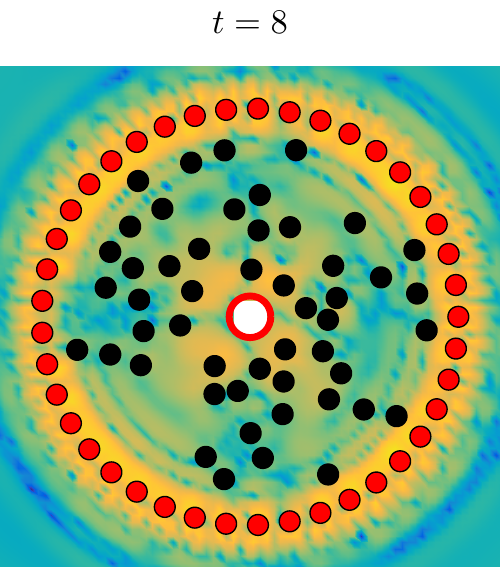}&
    \includegraphics{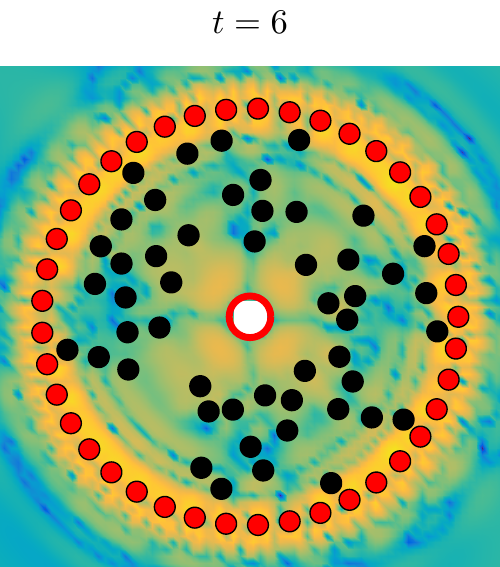}&
    \includegraphics{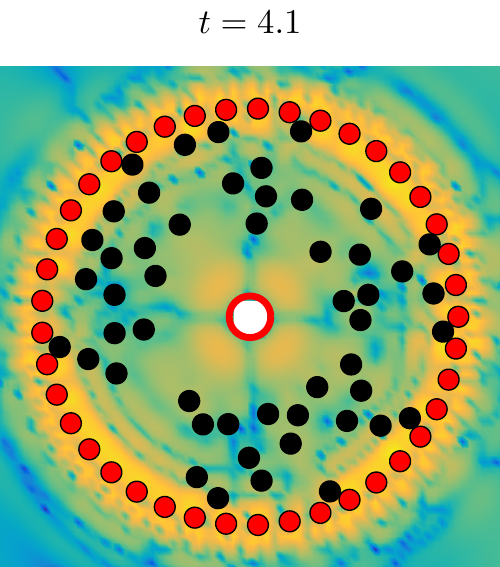}
  \end{tabular}
  \begin{center}
    \includegraphics{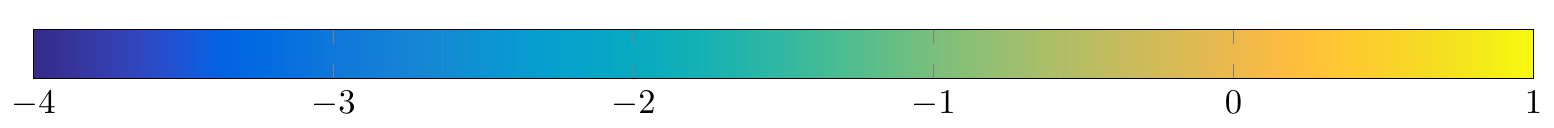}
  \end{center}
  \caption{\label{fig:macminn_stress} The shear strain rate
  $\log_{10}(|\sigma_{xy}|)$  of the suspension in
  Figure~\ref{fig:macminn}.  The formation of the petal-like patterns is
  also observed for the suspension of deformable
  bodies~\cite{MacMinn2015}.}
\end{figure}

In~\cite{MacMinn2015}, the shear strain rate is measured to better
characterize the flow.  In Figure~\ref{fig:macminn_stress}, we compute
and plot the shear strain rate for the simulation in
Figure~\ref{fig:macminn}.  A qualitative comparison of the numerical and
experimental results are in good agreement.  In particular, the
petal-like patterns in Figure~\ref{fig:macminn_stress} are also observed
in the experimental results.

\FloatBarrier
\subsection{Taylor-Couette Flow}
\begin{figure}[!h]
\begin{center}
\begin{tabular}{c c c c}
\includegraphics[width=3cm]{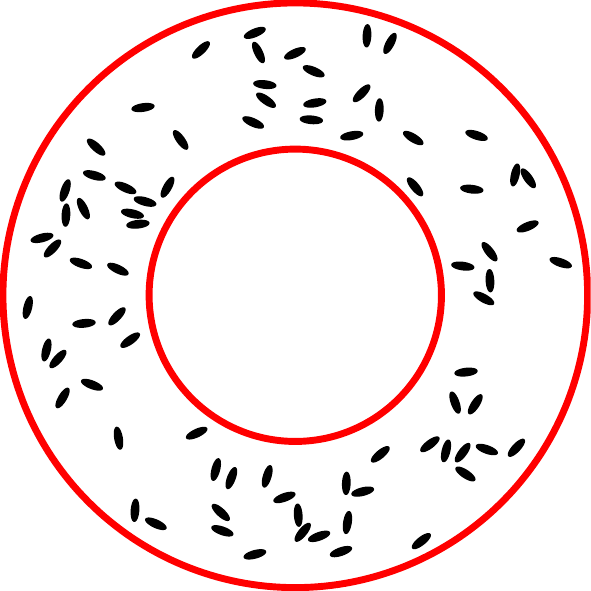} &
\includegraphics[width=3cm]{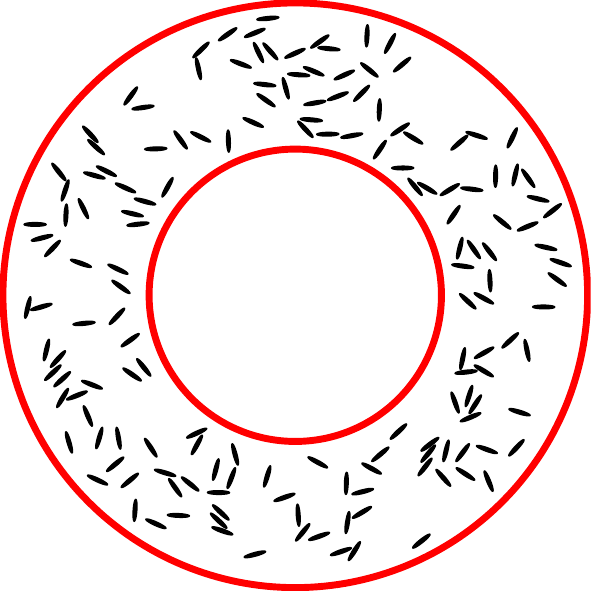} &
\includegraphics[width=3cm]{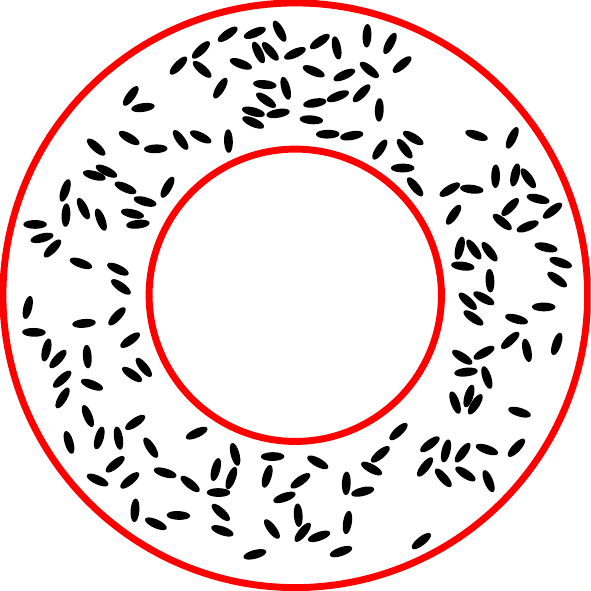} &
\includegraphics[width=3cm]{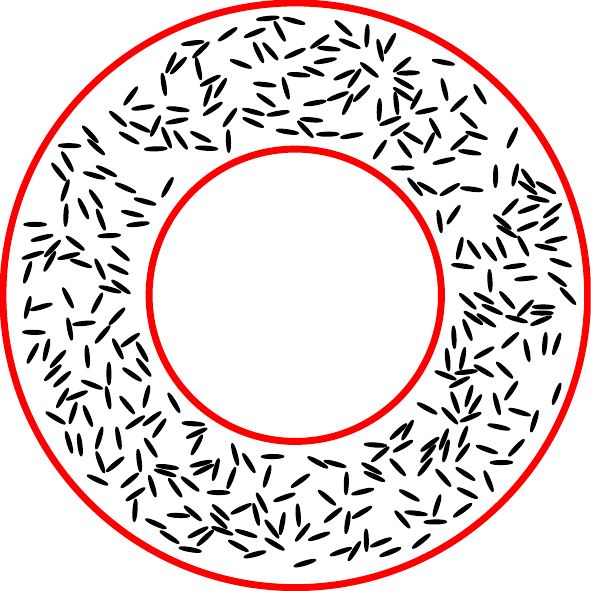}
\end{tabular}
\end{center}
\caption{Four initial configurations for Taylor-Couette flow with
varying volume fraction $\phi$ and aspect ratio $\lambda$. From left to
right: 1) $\phi=5\%$, $\lambda= 3$, 2) $\phi=5\%$, $\lambda=6$, 3)
$\phi=10\%$, $\lambda=3$, 4) $\phi=10\%$,
$\lambda=6$.}\label{fig:couette_setup}
\end{figure}

In many industrial applications, for example pulp and paper
manufacturing, suspensions of rigid elongated fibers are encountered.
Motivated by these suspensions we investigate rheological and
statistical properties of confined suspensions.  We consider suspensions
of varying area fraction and body aspect ratio; specifically we will
look at 5, 10, and 15 percent area fractions and elliptical bodies of
aspect ratio, $\lambda$ of 1, 3, and 6.  In all the examples, $nu M
1/ell^2$, so all the suspensions are in the dilute regime.  The bodies
initial locations are random, but non-overlapping (Figure
\ref{fig:couette_setup}).  The flow is driven by rotating the outer
cylinder at a constant angular velocity while the inner cylinder remains
fixed. 

Before measuring any rheological properties, we complete one full
revolution of the outer wall so that the bodies are well-mixed and
approaching a statistical equilibrium.  We start by considering the
alignment of the bodies. The alignment is particularly insightful since
many industrial processes involve fibers suspended in a flow, and the
alignment affects the material properties~\cite{larsoncf}.  One way to
measure the alignment is the order parameter, $S$ defined as,
\begin{align*}
  S = \left\langle \frac{d \cos^2\tilde{\theta} - 1}{d - 1} \right\rangle,
\end{align*}
where $d$ is the dimension of the problem ($2$ in our case),
$\tilde{\theta}$ is the deviation from the expected angle, and $\langle
\cdot\rangle$ averages over all bodies.  If $S=1$, all bodies are
perfectly aligned with the shear direction, $S=0$ corresponds to a
random suspension (no alignment), and $S=-1$ means that all bodies are
perfectly aligned perpendicular to the shear direction.  In our
geometric setup, a body centered at $(x,y)$ has an expected angle of
$\tilde{\theta}\tan^{-1}(y/x) + \pi/2$, and the average alignment of the
bodies will be in the direction of the shear, which is also
perpendicular to the radial direction.

\begin{figure}[!h]
\begin{center}
\includegraphics{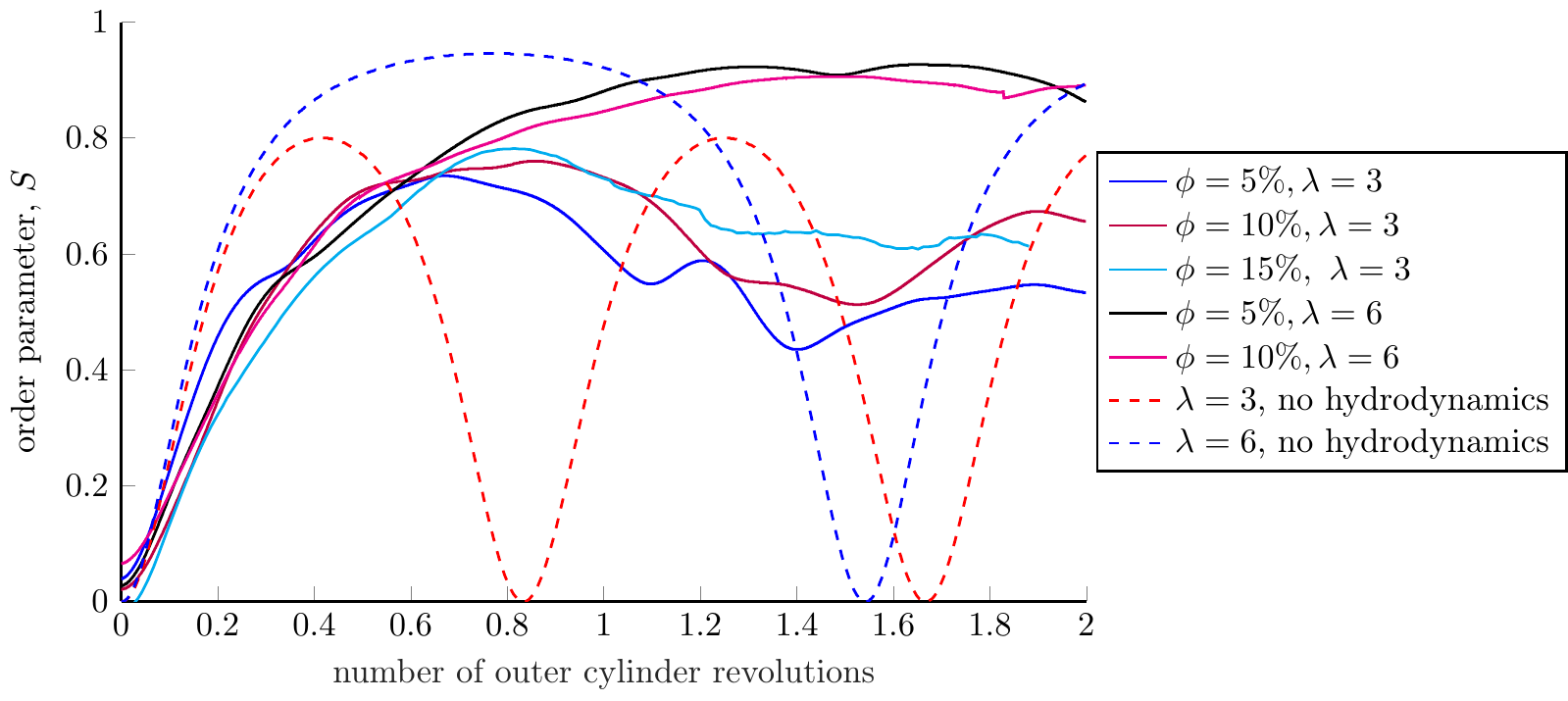}\\
\end{center}
\caption{\label{fig:angles} The order parameter of different fiber
concentrations and aspect ratios. We see that the 6:1 fibers align
better. The 6:1 fibers rotate rotate through the angle perpendicular to
the shear direction more quickly than the 3:1 fibers and thus spend more
time approximately aligned with the shear direction. The dashed lines
represent the order parameter for a suspension in an unbounded shear
flow with bodies that do not interact hydrodynamically. The red line
shows $\lambda=3$, while the blue line shows $\lambda=6$. }
\end{figure}

Since the initial condition is random, the initial configurations in
Figure~\ref{fig:couette_setup} have an order parameter $S\approx 0$. As
the outer cylinder rotates, we see in Figure~\ref{fig:angles} that $S$
increases quite quickly. The area fractions $\phi$ we consider have a
minor effect on $S$; however, the aspect ratio has a large effect.  In
particular, suspensions with slender bodies align much better with the
flow.

This matches the known dynamics of a single body in an unbounded shear
flow, where the body will align with the shear direction on
average. Bodies with a high aspect ratio rotate quickly when then they are
perpendicular to the shear direction and spend more time nearly aligned
with the shear direction. We compare our results to the time averaged
order parameter of a single elliptical body in an unbounded shear flow.
If the shear rate is $\dot{\gamma}$, the body rotates with period $\tau
= \pi/(2| \dot{\gamma}|)(\lambda + \lambda^{-1})$~\cite{Jeffery1922}
according to 
\begin{align*}
  \varphi(t) ~=~ \tan^{-1}\left(\frac{1}{\lambda}\tan\left(
    \frac{\lambda \dot{\gamma}t}{\lambda^2 + 1}\right)\right).
\end{align*}
The time average order parameter is then,
\begin{align*}
  \langle S\rangle ~=~ \frac{1}{\tau}\int_0^\tau\left( 
    2\cos^2(\varphi(t)) - 1\right)~\text{d}t ~=~ \frac{\lambda -1}{\lambda+1}.
\end{align*}
Independent of shear rate, for $\lambda= 3$ the theoretical $\langle
S\rangle$ is 1/2 and for $\lambda=6$ it is  5/7. Table \ref{tab:order}
shows the time and space averaged order parameter for the Couette
apparatus.We see that in all cases our computed time averaged order
parameter is higher than the theoretical single fiber case. This could
be due to the hydrodynamic interactions between the bodies, or the
effect of the solid walls.

In the absence of solid walls and hydrodynamic interactions between
bodies, a suspension will align and disalign. The period of the order
parameter in this case is the same as the rotational period for a single
fiber. In Figure \ref{fig:angles} the theoretical order parameter is
shown for a suspension of non-hydrodynamically interacting fibers in an
unbounded shear flow.  Hydrodynamic interactions prevent the suspension
from disaligning completely.

\begin{table}[!h]
\begin{center}
\begin{tabular}{c |c |c |c}
area fraction, $\phi$ & aspect ratio, $\lambda$ & computed $ \langle S \rangle$ & theoretical
$\langle S \rangle$ (single fiber)\\
\hline
5\%  & 3 & 0.52 & 0.50 \\
10\% & 3 & 0.60 & 0.50 \\
15\% & 3 & 0.65 & 0.50 \\
5\%  & 6 & 0.91 & 0.71 \\
10\% & 6 & 0.89 & 0.71
\end{tabular}
\end{center}
\caption{The time averaged order parameter during the second revolution of the Couette apparatus. The higher aspect ratio fibers align better on average. The alignment is in all cases higher than predicted for a single Jeffery orbit.
}\label{tab:order}
\end{table} 

\begin{figure}[!h]
\begin{center}
\includegraphics{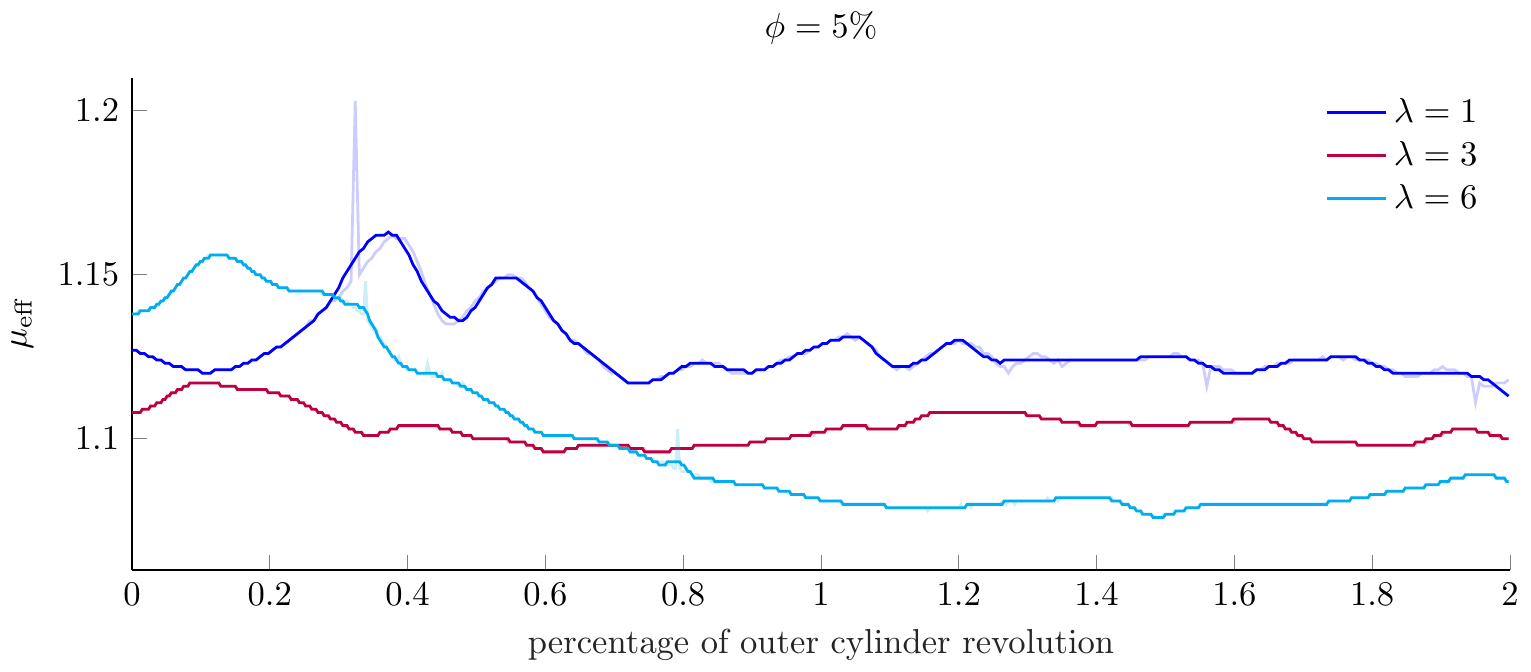}\\
\includegraphics{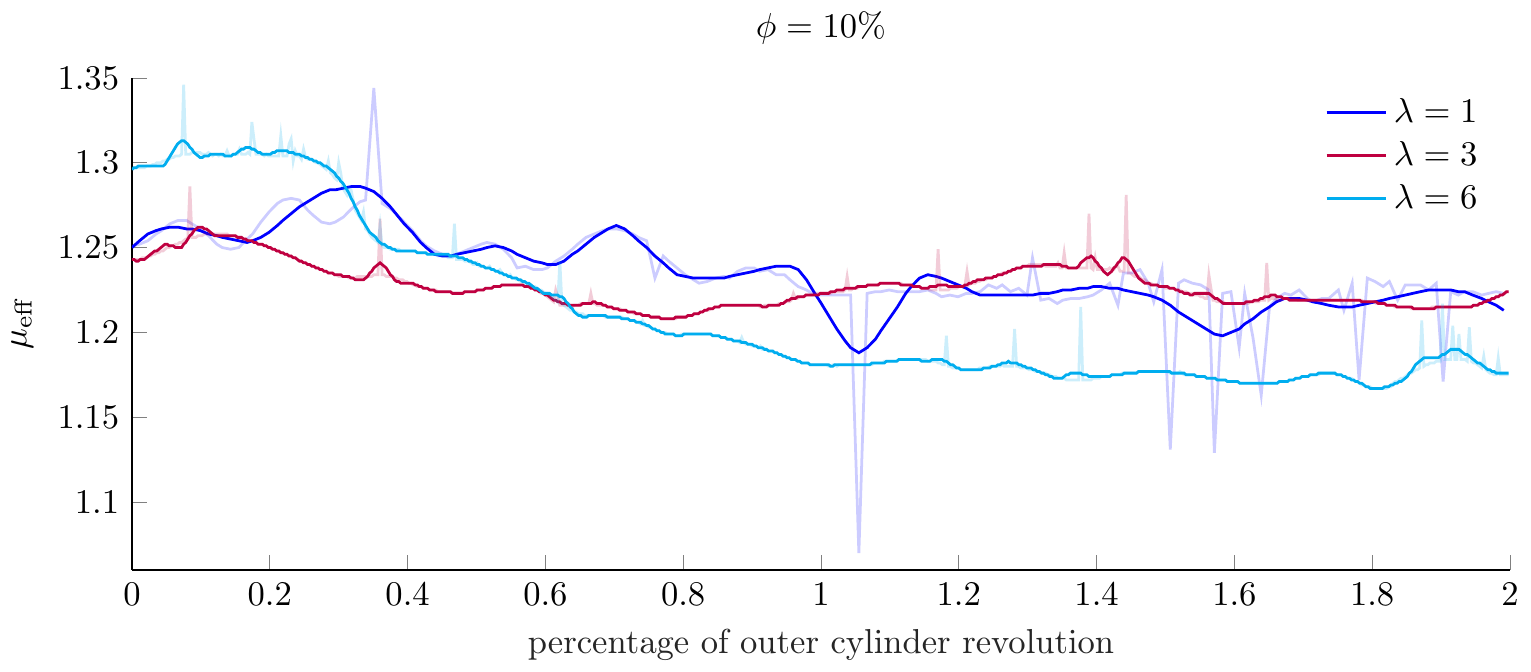}\\
\includegraphics{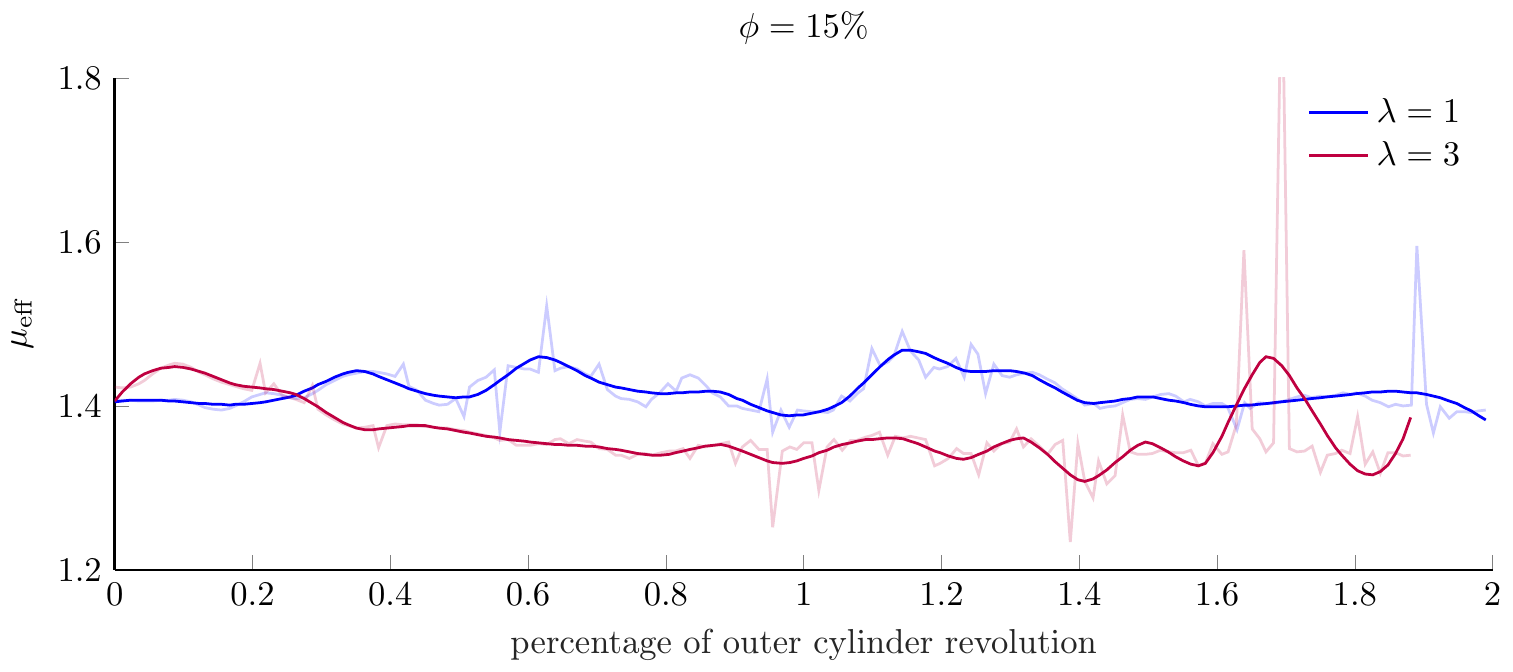}
\end{center}
\caption{Instantaneous bulk effective viscosity for various volume
fractions and body aspect ratios. The inner
cylinder is fixed, while the outer one rotates at a constant angular velocity. The transparent lines represent the raw data, while the solid lines have been smoothed using a multiscale local polynomial transform. }\label{fig:torque}
\end{figure}

Another quantity of interest in rheology is the effective viscosity of a
suspension. The shear viscosity $\mu$ relates the bulk shear stress
$\sigma_{xy}$ of a Newtonian fluid to the bulk shear rate
$\dot{\gamma}$, 
\begin{align*}
  \sigma_{xy} = \mu\dot{\gamma}.
\end{align*}
Adding bodies increases the bulk shear stress of a suspension.  The
proportionality constant relating the increased $\sigma_{xy}$ to the
shear rate is the {\em apparent viscosity}, and the ratio between the
apparent viscosity and the bulk viscosity is the effective viscosity
$\mu_{\text{eff}}$.  Experimentally, the bulk shear stress is often
computed by measuring the torque on the inner cylinder~\cite{Koos2012}.
Numerically, this is simply the strength of the rotlet centered in the
inner cylinder.  By computing the ratio of the torque on the inner
cylinder with bodies to the torque without bodies we determine the
effective viscosity of a suspension.  Figure~\ref{fig:torque} shows the
effective viscosity increases with $\phi$, but is generally lower for
bodies with aspect ratio $\lambda=6$. This is because higher aspect ratio
bodies align themselves better, and thus contribute less to the bulk
shear stress. The spikes in \ref{fig:torque} occur when a repulsion
force is added to the system. Similar spikes were observed numerically
in~\cite{Lu2017}. To make the results more clear we have used a multiscale local polynomial transform to smooth the data shown in Figure \ref{fig:torque}.

Finally, instead of computing the instantaneous effective viscosity,
experimenters are interested in the time averaged effective viscosity of
a suspension. In Table~\ref{tab:viscosity}, we report the average
instantaneous effective viscosity over the second revolution of the
outer cylinder
\begin{table}[!h]
\begin{center}
\begin{tabular}{c| c | c |  c}
$\lambda$ & 5\% area fraction  & 10\% area fraction  & 15\% area fraction\\
\hline
1 & 1.12 & 1.22 & 1.42 \\
3 & 1.10 & 1.23 & 1.36 \\
6 & 1.08 & 1.18 & 
\end{tabular}
\end{center}
\caption{Time averaged effective viscosity for various area fractions
and aspect ratios. The time average is done between the first and second
revolutions of the outer cylinder. As $\phi$ increases the effective
viscosity increases as expected. In general higher aspect ratio bodies
increase the viscosity less. }\label{tab:viscosity}
\end{table}

\FloatBarrier
\section{Conclusions\label{s:conclusions}}
We have developed a stable time stepping method for
simulating dense suspensions of rigid bodies without introducing
stiffness or requiring a large minimum separation distance.  A contact
algorithm is used that is guaranteed to avoid overlap.  Using this new
time stepper we are able to simulate concentrated suspensions with a
large time step and a small minimum separation distance.  We verify that
our new globally implicit time stepping method is able to take
larger stable time step sizes when compared to a locally implicit
time stepping method.

The contact algorithm is not reversible, and we examine the net effect
on the reversibility of the flow.  We compare the initial and final
configurations of two suspensions where the flow is reversed halfway
through the simulation.  The simulations demonstrate that the error in
the reversibility is the expected sum of the time stepping error and the
error introduced by the contact algorithm.

We use the new time stepping method to compute rheological properties of
suspensions including order parameters, effective viscosities, and
strain rates.  These rheological properties are compared with known
order parameters for of Jeffery orbits and the strain rate observed in
an experiment for deformable bodies.

There are several outstanding issues with this method. The most obvious
is that the suspensions are two-dimensional.  All the individual
algorithms we have developed in three dimensions.  However, to simulate
similar results with respect to volume fraction, minimum separation
distances, and aspect ratios, the spatial resolution will have to be
reduced, and this will still require parallel algorithms and
low-resolution algorithms~\cite{Kabacogulu2017}.  

A recent publication considers three-dimensional simulations of rigid
bodies~\cite{cor-gre-rac-vee2017}, but uses a different contact
algorithm than the STIV we used.  It is unclear that this contact
algorithm will be able to resolve dense suspensions where there are many
contacts.  As an alternative to the STIV, the contact can be measured
only at the final configuration, rather than as a space-time volume.
However, this collision metric indicates a contact-free configuration
when two particles pass through one another.  Therefore, we prefer to
use the STIV contact algorithm, which has been developed in three
dimensions~\cite{Harmon2011}, but it has not yet been applied to
suspensions.

Even with the development of our new time stepping method, we observe
difficult cases.  First, a good contact model between mobile bodies and
fixed rigid walls needs to be developed.  Second, the LCP solver can
converge slowly or not at all to a contact-free suspension.  While the
adaptive time stepping method helps alleviate some of these instances,
it is clear that additional algorithms are needed. Moreover, a more
reliable algorithm for choosing an optimal time step size would help
avoid instances where the time step size is reduced many times before a
time step is finally accepted.

\section*{References}


\end{document}